\documentclass[11pt]{article}

\usepackage{color}
\usepackage{latexsym}
\usepackage{amsmath,amsfonts,amssymb,theorem,euscript,array,enumerate,mathrsfs}

\usepackage[english,frenchb]{babel}
\usepackage[T1]{fontenc}

\newtheorem{Theorem}{Theorem}[section]
\newtheorem{Definition}{Definition}[section]
\newtheorem{Proposition}{Proposition}[section]

\newtheorem{Lemma}{Lemma}[section]
\newtheorem{Corollary}{Corollary}[section]
\newtheorem{Remark}{Remark}[section]
\newtheorem{Example}{Example}[section]
\newtheorem{Stasss}{Standing Assumption}[section]
\numberwithin{equation}{section}

\renewcommand{\abstractname}{le nom que tu souhaites afficher}

\def\esssup_#1{\underset{#1}{\mathrm{ess\,sup\, }}}
\def\essinf_#1{\underset{#1}{\mathrm{ess\,inf\, }}}

\def \trans{^{\scriptscriptstyle{\intercal}}}

\def \N{\mathbb{N}}
\def \R{\mathbb{R}}

\def \E{\mathbb{E}}
\def \F{\mathbb{F}}

\def \P{\mathbb{P}}

\def \S{\mathbb{S}}

\def \Fc{{\cal F}}
\def \Gc{{\cal G}}

\def \Lc{{\cal L}}

\def \eps{\varepsilon}

\def \ep{\hbox{ }\hfill$\Box$}

\def \epR{\hbox{ }\hfill$\lozenge$}

\def\reff#1{{\rm(\ref{#1})}}

\def\beqs{\begin{eqnarray*}}
\def\enqs{\end{eqnarray*}}
\def\beq{\begin{eqnarray}}
\def\enq{\end{eqnarray}}

\allowdisplaybreaks

\begin{document}

\title{BSDEs with diffusion constraint\\ and viscous Hamilton-Jacobi equations \\ with unbounded data
\thanks{We would like to thank the referees and AE for their suggestions which help us to improve the paper.}
}

\author{Andrea COSSO\thanks{Laboratoire de Probabilit\'es et Mod\`eles Al\'eatoires, CNRS, UMR 7599, Universit\'e Paris Diderot, \texttt{cosso@math.univ-paris-diderot.fr}}
\and
Huy\^{e}n PHAM\thanks{Laboratoire de Probabilit\'es et Mod\`eles Al\'eatoires, CNRS, UMR 7599, Universit\'e Paris Diderot, and CREST-ENSAE, \texttt{pham@math.univ-paris-diderot.fr}}
\and
Hao XING\thanks{Department of Statistics,
            London School of Economics and Political Science,
\texttt{h.xing@lse.ac.uk}}
}

\maketitle


\begin{abstract}
Nous donnons une repr\'esentation stochastique pour une classe g\'en\'erale d'\'equations d'Hamilton-Jacobi (HJ) visqueuses, convexes et super-nonlin\'eaires,  au moyen d'\'equations différentielles stochastiques rétrogrades (EDSR) avec contraintes sur la partie martingale. Nous comparons nos résultats avec la représentation classique en termes d'EDSR (super)quadratiques, et montrons notamment que l'existence d'une solution de viscosité à l'équation visqueuse de HJ peut 
être obtenue sous des conditions   de croi\-ssance plus générales, incluant des coefficients et une donnée  terminale non bornées.
\end{abstract}

\selectlanguage{english}


\begin{abstract}
We provide a stochastic representation for a general class of viscous Hamilton-Jacobi  (HJ)  equations, which has convex and superlinear nonlinearity in its gradient term, via a type of backward stochastic differential equation (BSDE) with constraint in the martingale part. We compare our result  with the classical representation in terms of  (super)quadratic BSDEs, and show in particular that existence of a viscosity solution to the viscous HJ equation can be obtained under more general  growth assumptions on the coefficients, including both unbounded diffusion coefficient and terminal data. 
\end{abstract}

\vspace{5mm}


\vspace{5mm}


\section{Introduction}

\setcounter{equation}{0}
\setcounter{Assumption}{0}
\setcounter{Theorem}{0}
\setcounter{Proposition}{0}
\setcounter{Corollary}{0}
\setcounter{Lemma}{0}
\setcounter{Definition}{0}
\setcounter{Remark}{0}
\setcounter{Example}{0}

Given $T\in(0,\infty)$ and $d\in \mathbb{N}\backslash\{0\}$, we consider the parabolic semilinear partial differential equation (PDE) of the form:
\beq \label{VHJ_0}
\left\{
\begin{array}{rcl}
- \dfrac{\partial u}{\partial t}(t,x)  - \frac{1}{2}\text{tr}\big[\sigma\sigma\trans(x)D_x^2 u(t,x)\big]- b(x)\cdot D_x u(t,x) & & \\
+ \;\; F\big(x,u(t,x),\varrho\trans(x)D_x u(t,x) \big) &= &0, \quad  (t,x)\in[0,T)\times\R^d, \\
u(T,x) &=& g(x),  \quad\, x\in\R^d.
\end{array}
\right.
\enq
Here $\sigma \trans$ and $\varrho\trans$  stand for the transpose of $\sigma$ and $\varrho$, $D_x u$ and $D^2_x u$ are the gradient and the Hessian of $u$, respectively, and $\cdot$ represents the inner product between two vectors in $\R^d$.
The functions $b\colon\R^d\rightarrow\R^d$, $\sigma,\varrho\colon\R^d\rightarrow\R^{d\times d}$,  and $F\colon\R^d\times\R\times\R^d\rightarrow\R$ are called coefficients of the equation, and $g\colon \R^d \rightarrow \R$ gives the terminal condition.

The focus of this paper is to provide a probabilistic representation for solutions to \eqref{VHJ_0}. This is a classical problem when $\varrho$ $=$ $\sigma$ and the standard approach  is to study a backward stochastic differential equation (BSDE) with the generator $F$ and a forward diffusion $X$, namely:
\beq \label{BSDEclassi}
\left\{
\begin{array}{ccl}
X_s^{t,x}  &=& x  +  \int_t^s b(X_r^{t,x}) dr +  \int_t^s \sigma(X_r^{t,x}) dW_r,  \;\;\;\;\;\;\;   (t,x) \in [0,T]\times\R^d, \; t \leq s \leq T, \\
Y_s^{t,x} &=& g(X_T^{t,x})  - \int_s^T F(X_r^{t,x},Y_r^{t,x},Z_r^{t,x}) dr - \int_s^T Z_r^{t,x} dW_r,
\end{array}
\right.
\enq
where $W$ is a $d$-dimensional Brownian motion. The forward backward stochastic differential equation \reff{BSDEclassi} is then formally connected to the  PDE \reff{VHJ_0} by the nonlinear Feynman-Kac formula: 
$Y_t^{t,x}$ $=$ $u(t,x)$.
Assuming  that $F$ is superlinear (which includes the quadratic and superquadratic case) and convex in its last argument (in which case, the PDE
\reff{VHJ_0} is referred to as the generalized \emph{viscous Hamilton-Jacobi equation}),
existence and uniqueness of a solution to \reff{BSDEclassi} is a quite complicated issue,  extensively  studied in the literature. This paper  provides   an alternative representation, which does not require the classical structural condition $\varrho$ $=$ $\sigma$, nor any non-degeneracy condition on $\varrho$ and $\sigma$.
In particular, our result can be applied to first-order equations, i.e. $\sigma$ $=$ $0$\footnote{When $\sigma=0$, the filtration $\mathbb{F}$ at the beginning of Section \ref{S:BSDE} below is generated by the Brownian motion
$B$ only.}.
Comparison to the literature will be presented in Remark \ref{remcomp} and in four examples below.
To state our main result, we introduce the  following assumptions which will be in force throughout the paper.

\begin{Stasss}\label{standing_assumption}
$\,$
\begin{enumerate}
\item[\textup{(i)}] There exists a constant $L_{b,\sigma,\varrho}$ such that
\[
|b(x) - b(x')| + \|\sigma(x) - \sigma(x')\| + \|\varrho(x) - \varrho(x')\|  \ \leq \ L_{b,\sigma,\varrho} |x - x'|, \;\; \text{for all } x, x' \in \R^d,
\]
where $\|A\|=\sqrt{\textup{tr}(AA\trans)}$ denotes the Frobenius norm of any matrix $A$.
\item[\textup{(ii)}] There exist  constants $M_\varrho$ $\geq$ $0$ and $p_\varrho$ $\in$ $[0,1]$ such that
\beqs
\|\varrho(x)\| &\leq& M_\varrho (1 + |x|^{p_\varrho}),  \quad \text{ for all } x \in \R^d.
\enqs
\item[\textup{(iii)}] There exists a constant $L_F$ such that
\[
 |F(x,y,z) - F(x,y',z)| \ \leq \ L_F |y - y'|, \quad \text{ for any } y,y'\in \R \text{ and } (x,z)\in \R^d \times \R^d.
\]
\item[\textup{(iv)}] The map $z\mapsto F(x,y,z)$ is convex. There exist constants $m_F, M_F\geq 0$, $q \geq p> 1$, $p_F\geq 0$, and $q_F\in [0, (1-p_\varrho)q/(q-1))$ if $p_\varrho$ $<$ $1$, or $q_F$ $=$ $0$ if $p_\varrho$ $=$ $1$,  such that
\[
- m_F \big(1 + |x|^{p_F} - \tfrac{1}{p}|z|^p\big) \ \leq \ F(x,0,z) \ \leq \ M_F \big(1 + |x|^{q_F} + \tfrac{1}{q}|z|^q\big),
\]
for any $(x,z)\in \R^d \times \R^d$.
\item[\textup{(v)}] Define the convex conjugate (or Fenchel-Legendre transform) of $F$ as
    \begin{equation}\label{def:f}
     f(x,a,y) = -\inf_{z\in \R^d}[a\cdot z + F(x,y,z)], \quad \text{ for all } (x,a,y) \in \R^d\times\R^d\times\R.
    \end{equation}
    The map $x\mapsto f(x,a,y)$ is continuous.
\item[\textup{(vi)}] The function $g$ is continuous. There exist constants $m_g, M_g, q_g\geq 0$,  and $p_g\in [0, (1-p_\varrho)q/(q-1))$  if $p_\varrho$ $<$ $1$, or $p_g$ $=$ $0$ if
$p_\varrho$ $=$ $1$, such that
\[
- m_g \big(1 + |x|^{p_g}\big) \ \leq \ g(x) \ \leq \ M_g \big(1 + |x|^{q_g}\big), \quad \text{ for all } x\in\R^d.
\]
\end{enumerate}
\end{Stasss}

 \vspace{1mm}

Let us now comment the above assumptions and their connection with related literature.

\begin{Remark} \label{remcomp}
{\rm {\bf 1.}  The condition $q$ $\geq$ $p$ $>$ $1$ in Assumption \ref{standing_assumption} (iv) means that $F$ has  a superlinear growth in $z$.
We also allow $F$ to have some polynomial growth in $x$, and  we distinguish the growth coefficients $p_F$ and $q_F$  for the lower and upper bounds. Indeed, notice that no condition is required on $p_F$, while we impose some upper  bound on $q_F$ depending  on $q$ and the (sub)linear growth coefficient $p_\varrho$
in Assumption \ref{standing_assumption} (ii).  Observe that this upper bound $\bar q_F$ $=$ $(1-p_\varrho)q/(q-1)$ is decreasing with $p_\varrho$ and $q$, with a limiting value equal to infinity when $q$ goes to $1$, and equal to $1-p_\varrho$ when $q$ goes to infinity; meanwhile $1-p_\varrho$ shrinks to zero (i.e. $F$ is upper-bounded in $x$) when $p_\varrho$ $=$ $1$ (i.e. $\varrho$ satisfies a linear growth condition).
Similarly, the terminal function $g$ is enabled to  satisfy a polynomial growth condition with the same constraint only on the power $p_g$ of the lower bound. These one-sided growth constraints on $F$ and $g$ are important for applications (see Example \ref{exa:utility} below) and are also sharp (see Example \ref{exa:Da Lio-Ley} below).

When $q=2$, $F$ has at most quadratic growth in $z$. A stochastic representation of \eqref{VHJ_0} with $\rho =\sigma$ is given by the \emph{quadratic} BSDE, which has been studied extensively, see \cite{kobylanski, briand_hu06, briand_hu08, delbaen_hu_richou11,richou12, barrieu_elkaroui, cheridito_nam,  bahlali_eddahbi_ouknine14}, amongst others.  For example, in the Markovian case of \cite{kobylanski}, $F$ and $g$  are assumed to be bounded in
$x$, i.e. $q_F$ $=$ $p_F$  $=$ $q_g$ $=$ $p_g$  $=$ $0$, moreover $\varrho$ $=$ $\sigma$ satisfies a linear growth condition, i.e. $p_\varrho$ $=$ $1$.
This case is covered by our assumptions (actually, we only need to assume $q_F$ $=$ $p_g$ $=$ $0$ when
$p_\varrho$ $=$ $1$,  but  no condition is required on $p_F$, $q_g$).  In the Markovian setting of \cite{briand_hu08}, \cite{delbaen_hu_richou11} or \cite{richou12},  $\varrho$ $=$ $\sigma$ is assumed to be a bounded function, i.e. $p_\varrho$ $=$ $0$. Notice that we do not assume any non-degeneracy condition on $\sigma$ so that the case of time dependent coefficients as in \cite{delbaen_hu_richou11} or \cite{richou12}, can be embedded in our framework by extending the spatial  variables from $x$ to $(t,x)$.  When $p>2$ in Assumption \ref{standing_assumption} (iv), $F$ has
\emph{super-quadratic} growth in $z$. This case has been studied in \cite{Ben-Artzi_Souplet_Weissler, gilding_guedda_kersner03, gladkov_guedda_kersner08, delbaen_hu_bao11, richou12, masiero_richou13}, which will be compared with our results in Example \ref{exa:super-QBSDE} below. Together with three other examples, we shall illustrate in further detail the scope of Assumption 1.1 and compare our conditions
on $q_F, p_g$  to existing results from both analytic and probabilistic aspects.

\noindent {\bf 2.} The case where $z\mapsto F(x,y,z)$ is \emph{concave} can be deduced from the convex case. Indeed, set $\tilde F(x,y,z)=-F(x, -y,-z)$. Then $\tilde F$ is convex in $z$ and our results apply to equation \eqref{VHJ_0} with $\tilde F$ in place of $F$ and $-g$ in place of $g$.

\noindent {\bf 3.} Assumption \ref{standing_assumption} (v) is satisfied under the following two sufficient conditions:
\begin{itemize}
\item The map $x \mapsto F(x,y,z)$ is continuous, uniformly with respect to $y$ and $z$, i.e., for all $x\in\R^d$ and $\eps>0$, there exists $\delta=\delta(x,\eps)>0$ such that,
\[
|F(x,y,z) - F(x',y,z)| \ \leq \ \eps, \quad \text{for any }|x - x'|\leq\delta \text{ and } (y,z)\in\R\times \R^d.
\]
\item The map $(x,z)\mapsto F(x,y,z)$ is continuous for any $y\in \R$, and there exists an optimizer $z^*= z^*(x,a,y)$ for  \eqref{def:f} such that $z^*$ is continuous in $x$.
\epR
\end{itemize}
}
\end{Remark}

 \vspace{1mm}

\begin{Example}[Generalized deterministic KPZ equation]\label{exa:KPZ}
{\rm
 Consider $F(z) = -\lambda |z|^q$ for some constants $\lambda>0$ and $q>1$. The equation \eqref{VHJ_0}, with $\varrho$ $=$ $\sigma$ as the identity matrix, is referred to as the \emph{generalized deterministic KPZ equation} with the $q=2$ case introduced by Kardar, Parisi, and Zhang in connection with the study of growing surfaces. In this case, $\tilde{F}(z) = \lambda |z|^q$ satisfies Assumption \ref{standing_assumption}, in particular,  $p_\varrho$ $=$ $0$ in (ii) and  $p_F=q_F =0$ in (iv). This equation has been studied from the analytical point of view in, for example, \cite{Ben-Artzi_Souplet_Weissler}, \cite{gilding_guedda_kersner03}, and \cite{gladkov_guedda_kersner08}. In particular, concerning Assumption \ref{standing_assumption} (iv), in \cite[Theorem 2.6] {gladkov_guedda_kersner08} $\max\{p_g, q_g\}<q/(q-1)$ is assumed, while here we only assume: $p_g$ $<$ $q/(q-1)$, but no restriction on $q_g$.
 Moreover, when $\max\{p_g, q_g\}= q/(q-1)$, an example whose solution explodes in finite time is presented in \cite[Remark 4.6]{gladkov_guedda_kersner08}. Together with the uniqueness result in \cite[Theorem 3.1]{gladkov_guedda_kersner08}, this example shows that global existence of solutions satisfying $|u(t,x)|\leq C(1+|x|^{q/(q-1)})$ is in general not expected when $p_g=q/(q-1)$.
\epR
}
\end{Example}

\begin{Example}[BSDEs with superquadratic growth]\label{exa:super-QBSDE}
{\rm Motivated by the previous exam\-ple, Delbaen, Hu, and Bao \cite{delbaen_hu_bao11} studied BSDEs whose generator has superquadratic growth in the ``$Z$-component". A solution to this superquadratic BSDE  provides a probabilistic representation for the solution of \eqref{VHJ_0} with $\varrho$ $=$
$\sigma$ and $q>2$. In particular, given a forward process
\[
 dX_s = b(s, X_s) ds + \sigma dW_s,
\]
with a Brownian motion $W$, a continuous differentiable function $b$ with bounded derivative, and a constant matrix $\sigma$, consider the following BSDE
\[
 Y_s = g(X_T) - \int_s^T F(Z_u) du -  \int_s^T Z_u dW_u,
\]
where $g$ is bounded and continuous, $F$ is nonnegative, convex, and satisfies $F(0)=0$ and $\overline{\lim}_{|z|\rightarrow \infty} \frac{F(z)}{|z|^2}=\infty$. Proposition 4.4 in \cite{delbaen_hu_bao11} presents a solution to the previous superquadratic BSDE. This result has been extended in \cite{richou12} and \cite{masiero_richou13}, where $F$ can depend on $x$ and $y$, without convexity assumption on $z$, and $\sigma$ can be a deterministic function of time. In these cases, $p_\varrho=0$, \cite[Assumptions (B.1) and (TC.1)]{masiero_richou13} implies that $\max\{p_F, q_F\}<q/(q-1)$ and $\max\{p_g, q_g\}<q/(q-1)$. When $\sigma$ depends on
$x$, only existence for small time is available, see \cite[Proposition 3.1]{richou12}. Compared to aforementioned works, our main result (see Theorem \ref{T:Feynman-Kac} below) provides an alternative representation for a global solution of \eqref{VHJ_0} when $F$ is convex in $z$, $\varrho$, not necessarily equals to $\sigma$, could depend on $x$ and have (sub)linear growth. Moreover no restrictions are imposed on $p_F$ and $q_g$. In particular, when $\varrho$ $=$ $\sigma$ is bounded, i.e. $p_{\varrho}=0$, our result assumes only $q_F, p_g<1$.
In fact, the asymmetry between upper and lower bounds of $F$ and $g$ is rather natural from BSDE point of view. Consider the BSDE
\[
 Y_s = g(W_T) -\int_s^T \tfrac12 |Z_u|^2 du -\int_s^T Z_u dW_u,
\]
whose solution is explicitly given by $Y_t = -\log \mathbb{E}[\exp(-g(W_T))]$. The previous expectation is well defined when $g$ is bounded from below by a sub-quadratic growth function, however no growth constraint on the upper bound of $g$ is needed. This asymmetry in assumptions also appeared in \cite{DHR} recently.
\epR
}
\end{Example}

\begin{Example}[Utility maximization]\label{exa:utility} {\rm Our  framework allows us  to incorporate the two follo\-wing financial applications.

\noindent (i)  {\it Portfolio optimization.} Consider a factor model of a financial market with a risk free asset $S^0$ and risky assets $S= (S^1, \cdots, S^n)$ with dynamics
 \beq \label{dynS}
  dS^0_t  \; = \;  S^0_t r dt & \text{and} &
  dS_t \; =  \; \text{diag}(S_t) [(r 1_n + \mu(X_t)) dt + a(X_t) \, dB_t],
 \enq
where $r\in \R$, $\mu$ and $a$ are measurable functions on $\R^d$, valued respectively on $\R^n$ and $\R^{n\times n}$, such that $a a\trans$ is invertible,
$\text{diag}(S)$ is a diagonal matrix with elements of $S$ on the diagonal, $1_n$ is a $n$-dimensional vector with every entry $1$, and $B$ is a $n$-dimensional Brownian motion.
We denote by $\lambda$ $=$ $a\trans(aa\trans)^{-1}\mu$,  the so-called Sharpe ratio.
In the dynamics \reff{dynS}, the factor $X$ is a  $d$-dimensional process governed by
 \beq \label{dynXfactor}
dX_t & = &  \beta(X_t) dt + \sigma(X_t) dW_t
 \enq
where $\beta$ and $\sigma$ are Lipschitz functions on $\R^d$, valued respectively on $\R^d$ and $\R^{d\times d}$,
and $W$ is a $d$-dimensional Brownian motion. The correlation between $B$ and $W$ is given by $d\langle W, B\rangle_t = \rho \,dt$ for some $\rho\in \R^{n\times d}$.

 An agent with power utility $U(w)= w^\gamma/\gamma$ for $\gamma<1$, $\gamma$ $\neq$ $0$ invests in this market in a self-financing way in order to maximize her expected utility of terminal wealth  at an investment horizon $T$. Let $v(t, w, x)$ be the value function of the investor and define the reduced value function $u$ via $v(t,w,x)= (w^\gamma/\gamma) e^{u(t,x)}$. Then the following equation, satisfied by $u$, is of the same type as \eqref{VHJ_0} with $\varrho$ $=$ $\sigma$ (see e.g. \cite[Equation (2.14)]{Nagai-06}):
 \beqs
\left\{
\begin{array}{rcl}
- \dfrac{\partial u}{\partial t}  - \frac{1}{2}\text{tr}\big[\sigma\sigma\trans D_x^2 u  \big]- b \cdot D_x u + F(x,\sigma\trans D_x u) & = & 0  \\
u(T,.) &=& 0,
\end{array}
\right.
\enqs
 where
 \begin{align*}
   b(x) \;&= \;  \beta(x) + \tfrac{\gamma}{1-\gamma} \sigma(x) \rho\trans\lambda(x), \qquad\qquad F(x, z) \; = \;  -\tfrac{1}{2} z M z\trans - h(x),\\
   M \; &= \;  \text{diag}(1_d) +  \tfrac{\gamma}{1-\gamma} \rho\trans\rho, \qquad\qquad\;\;\;\, h(x) \; = \;  \gamma \,r +  \tfrac{1}{2}\tfrac{\gamma}{1-\gamma} |\lambda(x)|^2.
 \end{align*}
 We assume that $b$ is Lipschitz. Moreover, note that $M$ is positive definite. Hence
 $\tilde{F}(x,z) = -F(x, -z)$ satisfies Assumption \ref{standing_assumption} (iv)  with $p$ $=$ $q$ $=$ $2$, and whenever $\lambda$ satisfies the condition:
 \beqs
 -m_F(1 + |x|^{p_F})  \; \leq \; \gamma |\lambda|^2 & \leq & M_F(1+|x|^{q_F}),
 \enqs
 for some nonnegative constants $m_F$, $M_F$, $p_F$, and $q_F$ $<$ $2(1-p_\varrho)$.  Hence, when $\gamma$ $>$ $0$, such condition  holds whenever $\lambda$ satisfies a strict sub-linear growth condition;
 while for $\gamma$ $<$ $0$ (the empirically relevant case in financial context), the previous condition holds once $\lambda$ satisfies a polynomial growth condition. In particular, the second scenario  includes the case where $a$ constant, and $\mu(x)$ (thus  $\lambda(x)$) is affine in $x$, as in the original Kim-Omberg model  where $X$ is an
Ornstein-Uhlenbeck process with $\sigma$ constant and $\beta(x)$ affine in $x$.

\vspace{1mm}

\noindent (ii) {\it Indifference pricing.}   We consider a financial model as in \reff{dynS}-\reff{dynXfactor}, where the process $X$  represents now  the level of  nontraded assets (e.g. volatility index, temperature), correlated with the traded assets of price $S$, for which the Sharpe ratio $\lambda$ is assumed to be bounded.  Given an European option
written on the  nontraded asset, with payoff $g(X_T)$ at maturity $T$, and  following the indifference pricing criterion (see e.g. \cite{musiela-zariphopoulou}),
we consider the problem of an agent with exponential utility $U(w)$ $=$ $-e^{-\gamma w}$,  $\gamma$ $>$ $0$, who invests in the traded assets $S$ up to $T$ where he has to
deliver the option $g(X_T)$.  Let $v(t,w,x)$ be the value function of the agent, and define the reduced value function $u$ via: $v(t,w,x)$ $=$ $U(w-u(t,x))$.  Then, $u$ satisfies an equation of type
\reff{VHJ_0} with $\varrho$ $=$ $\sigma$:
\beqs
\left\{
\begin{array}{rcl}
- \dfrac{\partial u}{\partial t}  - \frac{1}{2}\text{tr}\big[\sigma\sigma\trans D_x^2 u  \big]- b \cdot D_x u + F(x,\sigma\trans D_x u) & = & 0  \\
u(T,.) &=& g,
\end{array}
\right.
\enqs
 where
 \beqs
 b(x) \;= \;  \beta(x) +   \sigma(x) \rho\trans\lambda(x),  & &  F(x, z) \; = \;  - \tfrac{1}{2} z (\text{diag}(1_d) -  \rho\trans\rho) z\trans  +  \gamma \,r +  \tfrac{1}{2}   \frac{|\lambda(x)|^2}{\gamma}.
 \enqs
 Since $\lambda$ is assumed to be bounded,  $\tilde{F}(x,z) = -F(x, -z)$ satisfies Assumption \ref{standing_assumption} (iv) with $p$ $=$ $q$ $=$ $2$, $p_F$ $=$ $q_F$ $=$ $0$. Moreover,  Assumptions
 \ref{standing_assumption} (ii) and (vi) enable  us to consider unbounded diffusion coefficient $\sigma$ and unbounded payoff function $g$,
 for example when $X$  is governed by the ``shifted" CEV model with $\sigma(x)$ $=$ $\sigma_0 (\sigma_1+x)^{p_\varrho}$, for some positive constants $\sigma_0$,
 $\sigma_1$, and $p_\varrho$ $\in$ $[0,1)$ (the introduction of the positive constant $\sigma_1$ ensures that $\sigma$ is a Lipschitz function),   and $g$ satisfies the
 growth condition (recall Remark \ref{remcomp} {\bf 2.}):
 \beqs
 - M_g(1+ |x|^{q_g}) \; \leq \;  g(x) & \leq & m_g(1 + |x|^{p_g}),
 \enqs
 for some nonnegative constants $m_g$, $M_g$, $q_g$, and $0$ $\leq$ $p_g$ $<$ $2(1-p_\varrho)$.
\epR
 }
\end{Example}

\begin{Example}[A $\varrho\neq \sigma$ case and beyond]\label{exa:Da Lio-Ley}
{\rm The following linear-quadratic control problem was studied by Da Lio and Ley \cite[pp. 75]{DaLio-Ley-SICON}. Given $\mathbb{R}^{d\times d}$-valued deterministic functions $A, B, C$, and $D$, consider a linear stochastic differential equation
\[dX_s = [A(s) X_s + B(s) \alpha_s] ds + [C(s) X_s + D(s)] dW_s, \quad X_t=x.\]
Here $\alpha$ is taken from $\mathcal{A}$ which is the set of $\R^d$-valued predictable measure controls. Given another $\mathbb{R}^{d\times d}$-valued deterministic function $Q$, a constant
matrix $S\in \mathbb{R}^{d\times d}$, and a positive constant $R$, consider the following linear-quadratic (LQ) problem
\[
 V(t, x) = \inf_{\alpha\in \mathcal{A}} \mathbb{E}\Big[\int_t^T [X_s^{\trans} Q(s) X_s + R |\alpha_s|^2] ds + X_T^{\trans} S X_T\Big].
\]
By noting  that
\[
 -\inf_{a \in \mathbb{R}^d} \big[(D_x u)\trans B(t) a + R |a|^2\big] = \tfrac{1}{4R} |B(t)\trans D_x u|,
\]
the Hamilton-Jacobi-Bellman equation associated to this LQ problem is written as
\begin{equation} \nonumber
\left\{
\begin{array}{rcl}
- \frac{\partial u}{\partial t}  - \tfrac{1}{2}\text{tr}\big[a(t,x)D_x^2 u\big]- A(t) x\cdot D_x u -x\trans Q(t)x  + \tfrac{1}{4R} |B(t)\trans D_x u|^2 &=&  0, \;  \mbox{ on } [0,T)\times\R^d, \\
u(T,x) &=& x\trans S x,  \quad\, x\in\R^d,
\end{array}
\right.
\end{equation}
where $a(t,x) = (C(t) x + D(t)) (C(t) x + D(t))\trans$. In this case, $\sigma(t,x) = C(t) x + D(t) \neq B(t) = \varrho(t,x)$.

More generally, equation \eqref{VHJ_0} was studied in \cite{dalioley11} via viscosity solution techniques. Comparing to Assumption \ref{standing_assumption}, \cite{dalioley11} restricts to bounded $\varrho$, i.e., $p_{\varrho} =0$, and $\max\{p_F, q_F\}\leq q/(q-1)$. When $\max\{p_g, q_g\}<q/(q-1)$, the existence of a viscosity solution to \eqref{VHJ_0} is established in \cite[Theorem 2.2]{dalioley11}. 
}
\epR
\end{Example}

\vspace{1mm}

Before presenting our main result, let us first present an equivalent formulation of \eqref{VHJ_0} in terms of an Hamilton-Jacobi-Bellman (HJB) equation. Recall the convex conjugate function $f$ in \eqref{def:f}. By convexity of $F$ in $z$, we then  have the following duality relationship
\begin{equation*}
\label{F=-inf_f}
F(x,y,z) \ = \ - \inf_{a\in\R^d} \big[a\cdot z + f(x,a,y)\big], \quad \text{for all }\,(x,y,z)\in\R^d\times \R\times\R^d.
\end{equation*}
Therefore equation \eqref{VHJ_0} can be rewritten as the following HJB equation:
\begin{equation}
\label{VHJ}
\begin{cases}
\displaystyle - \frac{\partial u}{\partial t}(t,x) - \inf_{a\in\R^d}\big[\Lc^a u(t,x) + f(x,a,u(t,x))\big] \ = \ 0, &\quad (t,x)\in[0,T)\times\R^d, \\
u(T,x) \ = \ g(x), &\quad x\in\R^d,
\end{cases}
\end{equation}
where
\begin{equation}\label{E:La}
\Lc^a u(t,x) \ = \ (b(x) + \varrho(x) a) \cdot D_x u(t,x) + \frac{1}{2}\text{tr}\big[\sigma\sigma\trans(x)D_x^2 u(t,x)\big].
\end{equation}

It is standard to relate equation \eqref{VHJ}, at least formally, to the following optimal stochastic control problem of a recursive type:
\begin{equation}\label{E:u_control}
u(t,x) \ = \ \inf_\alpha \E\bigg[\int_t^T f\big(X_s^{t,x,\alpha},\alpha_s,u(s,X_s^{t,x,\alpha})\big) ds + g(X_T^{t,x,\alpha})\bigg],
\end{equation}
where the infimum is taken over all $\R^d$-valued predictable processes $\alpha$ and the controlled diffusion $X^{t,x,\alpha}$ evolves according to the equation
\begin{equation}
\label{FSDE_Controlled}
X_s^{t,x,\alpha} \ = \ x + \int_t^s \big(b(X_r^{t,x,\alpha}) + \varrho(X^{t,x,\alpha}_r) \, \alpha_r\big) dr + \int_t^s \sigma(X_r^{t,x,\alpha}) dW_r, \qquad t\leq s\leq T,
\end{equation}
where $W$ is a $d$-dimensional Brownian motion.

Rather than studying the control problem \eqref{E:u_control} directly, following \cite{kharroubi_pham12}, we introduce a \emph{randomized control} formulation. For every $(t,x,a)\in[0,T]\times\R^d\times\R^d$, we consider the following forward system of stochastic differential equations:
\begin{equation} \label{FSDE_X}
\begin{cases}
X^{t,x,a}_s =  x + \int_t^s \big(b(X^{t,x,a}_r) + \varrho(X^{t,x,a}_r) \,I_r^{t,a}\big) dr + \int_t^s \sigma(X^{t,x,a}_r) dW_r,  \\
I^{t,a}_s  =  a + B_s - B_t,
\end{cases}
\end{equation}
where $B$ is a $d$-dimensional Brownian motion, independent of $W$. In the next section, we shall check, under  Assumption \ref{standing_assumption} (i) and (ii), that there exists a unique solution
$(X^{t,x,a}, I^{t,a})$. System \eqref{FSDE_X}  is the randomized version of the controlled dynamics \eqref{FSDE_Controlled}. More precisely, the randomization procedure is performed by introducing the independent Brownian motion $B$, which is the natural choice when the control process $\alpha$ takes values in the entire space $\R^d$ as in the present case. On the contrary, if the control process $\alpha$ is $A$-valued, for some compact subset $A$ of $\R^d$, then a natural randomization is carried out by means of an independent Poisson random measure $\mu$ on $[0,\infty)\times A$; see \cite{kharroubi_pham12}.

Now we introduce a stochastic representation 
for \eqref{VHJ} and \eqref{VHJ_0}, via the following \emph{BSDE with diffusion constraint}. Given $(t,x,a)\in[0,T]\times\R^d\times\R^d$, consider
\begin{equation}
\label{BSDE}
Y_s \ = \ g(X_T^{t,x,a}) + \int_s^T f(X_r^{t,x,a},I_r^{t,a},Y_r) dr - \big(K_T - K_s\big) - \int_s^T Z_r dW_r - \int_s^T V_r dB_r,
\end{equation}
for all $s\in[t,T]$, together with the constraint
\begin{equation}
\label{BSDE_Constraint}
V_s \ = \ 0, \qquad ds\otimes d\P\text{-a.e.}
\end{equation}
The presence of the constraint \eqref{BSDE_Constraint} forces the introduction of the nondecreasing process $K$. In Theorem \ref{T:Y} in Section \ref{S:BSDE} below, we construct and prove the existence of a unique maximal solution $(Y^{t,x,a}, Z^{t,x,a}, V^{t,x,a}, K^{t,x,a})$ to \eqref{BSDE}-\eqref{BSDE_Constraint}. This allows us to present the main result of this paper, whose proof is presented in Section \ref{S:Feynman-Kac}.

\begin{Theorem}
\label{T:Feynman-Kac}
For all $(t,x)\in[0,T]\times\R^d$ and $a,a'\in\R^d$,
\[
Y^{t,x,a}_t \ = \ Y^{t,x,a'}_t.
\]
Define
\begin{equation}
\label{Feynman-Kac}
u(t,x) \ := \ Y_t^{t,x,a}, \qquad \text{for all }\,(t,x)\in[0,T]\times\R^d.
\end{equation}
Then there exists a constant $C$ such that
\begin{equation}\label{growth u}
 |u(t,x)|\leq C(1+|x|^{p_F \vee q_F\vee p_g \vee q_g}), \quad \text{ for all } (t,x) \in [0,T]\times \mathbb{R}^d,
\end{equation}
and $u$ is a viscosity solution to equation \eqref{VHJ_0} $($or, equivalently, to \eqref{VHJ}$)$.
\end{Theorem}

\begin{Remark}
{\rm Theorem \ref{T:Feynman-Kac} focuses on global solutions of \eqref{VHJ_0}, i.e. $T$ is not necessarily small. Local solutions have also been studied in the literature. When $\varrho$ is bounded, i.e., $p_{\varrho}=0$, and $\max\{p_F, q_F, p_g, q_g\}\leq q/(q-1)$, \cite[Theorem 3.2]{dalioley11} establishes a solution $u$ satisfying 
\begin{equation}\label{growth}
 |u(t,x)| \leq C(1+ |x|^{q/(q-1)}), \quad \text{ for some } C,
\end{equation}
and equation \eqref{VHJ_0} on $[T-\tau, T]$ for some $\tau>0$. Similar local existence has also been obtained via probabilistic methods for (super)quadratic BSDEs in \cite[Proposition 4.2]{delbaen_hu_richou11} and \cite[Proposition 3.1]{richou12}. For global existence, the restriction on $p_g$ in Assumption \ref{standing_assumption} (vi) is actually sharp. When $p_{\varrho}=0$ and $p_g=q/(q-1)$, \cite[Exemple 3.3]{dalioley11} presents a deterministic control problem (i.e. $\sigma=0$) whose value function is the only possible viscosity solution of \eqref{VHJ_0} in the class of functions satisfying \eqref{growth}. However, under some parameter specification, the value function blows up when time is less than $T-\tau$ for some $\tau>0$. Therefore, similar to Example \ref{exa:KPZ}, when $p_g=(1-p_{\varrho})q/(q-1)$, global viscosity solutions of \eqref{VHJ_0} satisfying \eqref{growth} do not exist in general, and the growth constraint in Assumption \ref{standing_assumption} (vi) cannot be improved even when \eqref{VHJ_0} is a first order equation (i.e., $\sigma=0$).
}
\end{Remark}

\begin{Remark}
{\rm 
When $\varrho$ is bounded, $\max\{p_F , q_F, p_g, q_g\} \leq q/(q - 1)$, a comparison result for sub(super)solutions to \eqref{VHJ_0} in the class of functions satisfying the growth condition \eqref{growth u} was obtained in \cite[Theorem 3.1]{dalioley11}. As a result, $u$ given by \eqref{Feynman-Kac} is the unique solution to equation \eqref{VHJ_0} in the class of functions satisfying \eqref{growth u}. In this case, $u$ is also continuous, as a byproduct of the comparison result.
\epR
}
\end{Remark}

\begin{Remark}
{\rm  When there exists a smooth solution $u$  to the equation \eqref{VHJ_0} (or equivalently to the PDE \reff{VHJ}), one can check directly  the connection between PDE \eqref{VHJ_0} or \eqref{VHJ} and the BSDE with
diffusion constraint \reff{FSDE_X}-\reff{BSDE}-\reff{BSDE_Constraint}.  Indeed, by applying It\^o's formula to $u$ along the forward diffusion process governed by \reff{FSDE_X}, we see that
$Y_s$ $=$ $u(s,X_s^{t,x,a})$, $t\leq s\leq T$,  satisfies the relation \reff{BSDE} with
\beqs
K_s & = &  \int_t^s \Big[  \big(\frac{\partial u}{\partial t}  + \Lc^{I_r^{t,a}} u\big)(r,X_r^{t,x,a}) + f(X_r^{t,x,a},I_r^{t,a},Y_s) \Big] dr \\
Z_s &=&  \sigma\trans(X_r^{t,x,a})D_x u(s,X_s^{t,x,a}), \;\;\; V_s \; = \; 0, \;\;\; t \leq s \leq T.
\enqs
Since $u$ is a solution to the PDE \reff{VHJ}, the term inside the bracket of $K$ is nonnegative, which implies that $(K_s)_{t\leq s\leq T}$
is a nondecreasing process starting from $K_t$ $=$ $0$, and therefore the quadruple $(Y,Z,V,K)$ is solution to
the BSDE with diffusion constraint \reff{FSDE_X}-\reff{BSDE}-\reff{BSDE_Constraint}. Actually, this holds true whenever  $u$ is a subsolution to the PDE \reff{VHJ}, and the fact that $u$ is a solution to this PDE
will imply that  $(Y,Z,V,K)$ is {\it the maximal solution}  to the BSDE with diffusion constraint \reff{FSDE_X}-\reff{BSDE}-\reff{BSDE_Constraint}. 
This intuition will be proved in the rest of the paper without assuming the existence of a smooth solution to \eqref{VHJ_0}.}
\epR
\end{Remark}

\section{BSDE with diffusion constraint}
\label{S:BSDE}

\setcounter{equation}{0}
\setcounter{Assumption}{0}
\setcounter{Theorem}{0}
\setcounter{Proposition}{0}
\setcounter{Corollary}{0}
\setcounter{Lemma}{0}
\setcounter{Definition}{0}
\setcounter{Remark}{0}
\setcounter{Example}{0}

Let us introduce our probabilistic setting. Consider a filtered probability space $(\Omega,\Fc_T,\F=(\Fc_s)_{0\leq s\leq T},\P)$, where $\F$ is the standard augmentation of the filtration generated by two $d$-dimensional independent Brownian motions $W=(W_s)_{0\leq s\leq T}$ and $B=(B_s)_{0\leq s\leq T}$. For $0\leq t\leq T$, we also consider the ``shifted" version $(\Omega, \Fc^t_T, \F^t=(\Fc^t_s)_{t\leq s\leq T}, \P)$, where $(W_s)_{0\leq s\leq T}$ and $(B_s)_{0\leq s\leq T}$ before are replaced by $(W_s-W_t)_{t\leq s\leq T}$ and $(B_s-B_t)_{t\leq s\leq T}$, respectively. We denote by $\E_s^t$ the conditional expectation under $\P$ given $\Fc_s^t$ for $0\leq t\leq s\leq T$, and see that $\E_t^t$ coincides with the ``ordinary" expectation $\E$.
On this ``shifted" probability space, we introduce the following spaces of stochastic processes.

\begin{itemize}
\item $\mathbb S^2(t,T)$: the set of real-valued c\`adl\`ag $\F^t$-adapted processes $Y$ $=$
$(Y_s)_{t\leq s\leq T}$ satisfying
\[
\|Y\|_{_{\mathbb S^2(t,T)}}^2 \ := \ \E\Big[ \sup_{t\leq s\leq T} |Y_s|^2 \Big] \ < \ \infty.
\]
\item $\mathbb H^2(t,T)$: the set of $\R^d$-valued $\F^t$-predictable processes
$Z=(Z_s)_{t\leq s\leq T}$ satisfying
\[
\|Z\|_{_{\mathbb H^2(t,T)}}^2 \ := \ \E\bigg[\int_t^T |Z_s|^2 ds\bigg] \ < \ \infty.
\]
\item $\mathbb K^2(t,T)$: the set of nondecreasing $\F^t$-predictable processes $K$ $=$ $(K_s)_{t\leq s\leq T}$ $\in$  $\mathbb S^2(t,T)$ with $K_t$ $=$ $0$. We have
\[
\|K\|_{_{\mathbb S^2(t,T)}}^2 \ = \ \E\big[K_T^2\big].
\]
\end{itemize}

\vspace{2mm}

This section focuses on the construction of a maximal solution to the forward backward SDE with diffusion constraint
\reff{FSDE_X}-\eqref{BSDE}-\eqref{BSDE_Constraint}.  We first check the existence of a unique solution to the randomized forward system \reff{FSDE_X},  and show some useful moment estimates. Throughout this section, Assumption \ref{standing_assumption} is in force.

\begin{Lemma} \label{Estimate_X}
For any  $(t,x,a)\in[0,T]\times\R^d\times\R^d$, there exists a unique (up to indistinguishability) adapted and continuous process $\{(X_s^{t,x,a},I_s^{t,a}), t\leq s\leq T\}$ solving  \reff{FSDE_X}. Moreover,  for every $m\geq1$, there exists a positive constant $C$, depending only on $m$, $T$, $L_{b,\sigma,\varrho}$,  $M_\varrho$, and $p_\varrho$ such that
\begin{equation} \label{estimXprho<1}
\E_s^t\Big[\sup_{r\in[s,T]}|X_r^{t,x,a}|^m  \Big] \ \leq \
C\bigg(1 + |X_s^{t,x,a}|^m +  \int_s^T \E_s^t\big[|I_r^{t,a}|^{\frac{m}{1-p_\varrho}}\big] dr\bigg), \;\;  \mbox{ if } \;\;  p_\varrho  <1,
\end{equation}
and, if $p_\varrho=1$,
\begin{equation}
\E_s^t\Big[\sup_{r\in[s,T]}|X_r^{t,x,a}|^m\Big] \ \leq \
C\bigg(1 + |X_s^{t,x,a}|^m  +  \sqrt{\int_s^T \E_s^t\big[|I_r^{t,a}|^{2m}] dr} \bigg) \sqrt{\E_s^t\big[e^{C \int_s^T |I_r^{t,a}| dr }\big]},  \label{estimXprho=1}
\end{equation}
for all $0\leq t\leq s\leq T$.
\end{Lemma}
{\bf Proof.} For every $(t,a)\in[0,T]\times\R^d$, there exists a unique process $I^{t,a}=(I_s^{t,a})_{t\leq s\leq T}$ satisfying the second equation in \reff{FSDE_X}. Since the coefficients $b(x) + \varrho(x)a$ and $\sigma(x)$ of the diffusion system \reff{FSDE_X} for $(X,I)$  are locally Lipschitz in $(x,a)$, it is well-known (see e.g. Exercise IX.2.10 in \cite{revuz_yor}) that for every $(t,x,a)\in[0,T]\times\R^d\times\R^d$ there exists an adapted process $X^{t,x,a}=(X_s^{t,x,a})_{t\leq s\leq T}$ such that, if $e:=\inf\{s\geq t\colon|X_s^{t,x,a}|=\infty\}$ with $\inf\emptyset=\infty$, then $X^{t,x,a}$ is the unique (up to indistinguishability) adapted and continuous process on $[t,e)\cap[t,T]$ that satisfies the first equation in \reff{FSDE_X} on $[t,e)\cap[t,T]$. It then remains to prove that the explosion time $e$ of $X^{t,x,a}$ satisfies: $\P(e=\infty)=1$.

For simplicity of notation, we denote $I=I^{t,a}$ and $X=X^{t,x,a}$. Consider $T_n=\inf\{s\geq t\colon|X_s|>n\}\wedge T$ and apply It\^o's formula to $|X|^m$, for $m$ $\geq$ $1$. From the dynamics \reff{FSDE_X} of $(X,I)$ and Young's inequality, we see, under  Assumption  \ref{standing_assumption} (i) and (ii), that  there exists a constant $C$ (which in the sequel may change from line to line) depending only on $m$, $T$, $L_{b,\sigma,\varrho}$, $M_\varrho$, and $p_\varrho$, such that for all $t\leq s \leq T_n$,
\beq
\sup_{r\in [s,T_n]} |X_r|^m & \leq & C\bigg( 1 + |X_s|^m + \int_s^{T_n} \big( |X_r|^m +  |I_r|^m + |X_r|^{m - 1 + p_\varrho} |I_r| \big) dr  \nonumber \\
& &  \hspace{2cm}  + \;  \sup_{r\in [s,T_n]} \bigg| \int_s^r |X_u|^{m-1} \sigma(X_u) dW_u \bigg|  \bigg).  \label{Xm}
\enq

\textbf{Case $p_\varrho$ $<$ $1$.} By applying Young's inequality to $|X_r|^{m-1+p_\varrho} |I_r|$,  taking conditional expectation on both sides of \reff{Xm}, and using the Burkholder-Davis-Gundy inequality,
we get
\beqs
\E_s^t \bigg[ \sup_{r\in [s,T_n]} |X_r|^m  \bigg] & \leq & C \bigg( 1 + |X_s|^m +  \int_s^T \E_s^t \big[ |I_r|^{\frac{m}{1-p_\varrho}} 1_{[s,T_n]}(r) \big] dr \\
& & + \; \int_s^T \E_s^t\big[ \sup_{u\in[s,r]} (|X_u|^m 1_{[s,T_n]}(u)) \big] dr \bigg),
\enqs
which shows, by Gronwall's lemma,
\beqs
\E_s^t\Big[\sup_{r\in[s,T_n]}|X_r|^m  \Big] & \leq &
C\bigg(1 + |X_s|^m +  \int_s^T \E_s^t\big[|I_r|^{\frac{m}{1-p_\varrho}} 1_{[s,T_n]}(r) \big] dr\bigg) \\
& \leq &
C\bigg(1 + |X_s|^m +  \int_s^T \E_s^t\big[|I_r|^{\frac{m}{1-p_\varrho}}\big] dr\bigg).
\enqs
Since $T_n\nearrow e\wedge T$ as $n$ goes to infinity, from Fatou's lemma we obtain
\beq \label{estimXprho<1_e}
\E_s^t\Big[\sup_{r\in[s,e\wedge T]}|X_r|^m  \Big] & \leq & C\bigg(1 + |X_s|^m +  \int_s^T \E_s^t\big[|I_r|^{\frac{m}{1-p_\varrho}}\big] dr\bigg).
\enq
In particular, taking $s=t$, we have $\P(e=\infty)=1$. Then, from \reff{estimXprho<1_e} we deduce the required estimate \reff{estimXprho<1}.

\vspace{1mm}

\textbf{Case $p_\varrho$ $=$ $1$.} Take the conditional expectation in
\reff{Xm}  with respect to the $\sigma$-algebra $\Gc_s^t:=\Fc_s^t\vee\sigma(I_r,\,t\leq r\leq T)$, and observe that $(W_s-W_t)_{t\leq s\leq T}$ is a Brownian motion with respect to $(\Gc_s^t)_{t\leq s\leq T}$ since $W$ and $I$ are independent. Therefore, we can still use the Burkholder-Davis-Gundy inequality, and obtain (we denote by $\E_s^{t,\Gc}$ the conditional expectation under $\P$ given $\Gc_s^t$)
\beqs
\E_s^{t,\Gc} \Big[ \sup_{r\in [s,T_n]} |X_r|^m  \Big] & \leq & C \bigg( 1 + |X_s|^m +  \int_s^T |I_r|^m \E_s^{t,\Gc}\big[1_{[s,T_n]}(r)\big] dr \\
& & + \; \int_s^T  \E_s^{t,\Gc} \big[ \sup_{u\in [s,r]} (|X_u|^m 1_{[s,T_n]}(u)) \big] \big( 1 + |I_r|\big) dr \bigg),
\enqs
which gives, by Gronwall's lemma and noting that $T_n\leq T$,
\beqs
\E_s^{t,\Gc} \Big[ \sup_{r\in [s,T_n]} |X_r|^m \Big] & \leq & C \bigg( 1 + |X_s|^m +  \int_s^T |I_r|^m dr \bigg) \exp\bigg( C \int_s^T (1 + |I_r|) dr \bigg) \\
&\leq & C \bigg( 1 + |X_s|^m +  \int_s^T |I_r|^m dr \bigg) e^{C \int_s^T |I_r| dr}.
\enqs
Afterwards, by taking conditional expectation with respect to $\Fc_s^t$, using the Cauchy-Schwarz inequality, and the fact that $\sqrt{a+b} \leq \sqrt{a} + \sqrt{b}$ for any $a,b\geq0$, we get
\beqs
\E_s^t\Big[\sup_{r\in[s,T_n]}|X_r|^m\Big] & \leq &
C\bigg(1 + |X_s|^m  +  \sqrt{\int_s^T \E_s^t\big[|I_r|^{2m}] dr} \bigg)\sqrt{\E_s^t\big[e^{C \int_s^T |I_r| dr}\big]}. \nonumber
\enqs
Recalling that $T_n\nearrow e\wedge T$ as $n$ goes to infinity, from Fatou's lemma we obtain
\beq
\E_s^t\Big[\sup_{r\in[s,e\wedge T]}|X_r|^m\Big] & \leq &
C\bigg(1 + |X_s|^m  +  \sqrt{\int_s^T \E_s^t\big[|I_r|^{2m}] dr} \bigg) \sqrt{\E_s^t\big[e^{C \int_s^T |I_r| dr}\big]},  \label{estimXprho=1_e}
\enq
where the second conditional expectation on the right-hand side is finite, since $I$ is a Brownian motion. Therefore, \eqref{estimXprho=1_e} with $s=t$ implies  $\P(e=\infty)=1$. Finally, from \reff{estimXprho=1_e} we get the required estimate \reff{estimXprho=1}.
\ep

\vspace{3mm}

Next, it is straight forward to check that Assumption \ref{standing_assumption} translate to the following properties on the generator $f$ of the BSDE \eqref{BSDE}.

\begin{Lemma}
\label{L:Fenchel-Legendre}
The map $a\mapsto f(x,y,a)$ is convex and satisfies
\begin{align}
|f(x,a,y) - f(x,a,y')| \ &\leq \ L_F |y - y'|, \label{Lipschitz_f} \\
- M_F \big(1 + |x|^{q_F} - \tfrac{1}{q' M_F^{q'}}|a|^{q'}\big) \ \leq \ f(x,a,0) \ &\leq \ m_F \big(1 + |x|^{p_F} + \tfrac{1}{p' m_F^{p'}}|a|^{p'}\big), \label{GrowthCond_f}
\end{align}
for all $x\in\R^d$, $y,y'\in\R$, and $a\in\R^d$. Here $p'=p/(p-1)$ and $q'=q/(q-1)$ are the conjugate exponents of, respectively, $p$ and $q$ in Assumption \ref{standing_assumption}.
\end{Lemma}

\vspace{1mm}

We may then  define the notion of maximal solution to the   BSDE with diffusion constraint  \eqref{BSDE}-\eqref{BSDE_Constraint}.

\begin{Definition}
For every $(t,x,a)\in[0,T]\times\R^d\times\R^d$, we say that $(Y^{t,x,a},Z^{t,x,a},V^{t,x,a},$ $K^{t,x,a})$ $\in$ $\mathbb S^2(t,T)\times\mathbb H^2(t,T)\times\mathbb H^2(t,T)\times\mathbb K^2(t,T)$ is a maximal solution to \eqref{BSDE}-\eqref{BSDE_Constraint} if it satisfies \eqref{BSDE}-\eqref{BSDE_Constraint}, and for any other solution $(\underline Y^{t,x,a},\underline Z^{t,x,a},\underline V^{t,x,a},\underline K^{t,x,a})\in\mathbb S^2(t,T)\times\mathbb H^2(t,T)\times\mathbb H^2(t,T)\times\mathbb K^2(t,T)$ we have $Y^{t,x,a}_s\geq\underline Y^{t,x,a}_s$, for any $s\in[t,T]$.
\end{Definition}

Such a maximal solution is constructed
using a penalization approach as in \cite[Theorem 2.1]{kharroubi_pham12}. However, rather than employing an independent Poisson random measure as in \cite{kharroubi_pham12}, our randomization here is carried out by means of an independent Brownian motion in \eqref{FSDE_X}.
Let us consider, for every $(t,x,a)\in[0,T]\times\R^d\times\R^d$ and $n\in\N$, the following \emph{penalized} BSDE:
\begin{align}
\label{BSDE_n}
Y_s^n \ &= \ g(X_T^{t,x,a}) + \int_s^T f(X_r^{t,x,a}, I_r^{t,a}, Y_r^n) dr - \big(K_T^n - K_s^n\big) - \int_s^T Z_r^n dW_r - \int_s^T V_r^n dB_r \notag \\
&= \ g(X_T^{t,x,a}) + \int_s^T f_n(X_r^{t,x,a},I_r^{t,a}, Y_r^n, V_r^n) dr - \int_s^T Z_r^n dW_r - \int_s^T V_r^n dB_r,
\end{align}
for all $s\in[t,T]$, where
\[
K_s^n \ = \ n \int_t^s |V_s^n| ds
\]
and the generator $f_n(x,a,y,v) \ = \ f(x,a,y) - n|v|$. Notice by \reff{Lipschitz_f} that this generator $f_n$ is Lipschitz in $(y,v)$, so that from standard result due to \cite{pardoux_peng90}, we know that for every $(t,x,a)\in[0,T]\times\R^d\times\R^d$ and $n\in\N$, there exists a unique solution
$(Y^{n,t,x,a},Z^{n,t,x,a},V^{n,t,x,a})\in\S^2(t,T)\times\mathbb H^2(t,T)\times\mathbb H^2(t,T)$ to BSDE \eqref{BSDE_n}.

Moreover, we have the following comparison results.
\begin{Lemma}
\label{L:Monotone}
For every $(t,x,a)\in[0,T]\times\R^d\times\R^d$, the following statements hold:
\begin{enumerate}
\item[\textup{(i)}] The sequence $(Y^{n,t,x,a})_n$ is nonincreasing, i.e., $Y^{n,t,x,a}\geq Y^{n+1,t,x,a}$, $n\in\N$.
\item[\textup{(ii)}] For any solution $(\bar Y^{t,x,a},\bar Z^{t,x,a},\bar V^{t,x,a},\bar K^{t,x,a})\in\S^2(t,T)\times\mathbb H^2(t,T)\times\mathbb H^2(t,T)\times\mathbb K^2(t,T)$ to \eqref{BSDE}-\eqref{BSDE_Constraint}, we have $Y^{n,t,x,a}\geq\bar Y^{t,x,a}$, $n\in\N$.
\end{enumerate}
\end{Lemma}
\textbf{Proof.} Since $f_n\geq f_{n+1}$, the first statement follows from a direct application of the comparison theorem for BSDEs. For the second statement, note that $f_n(X_r^{t,x,a},I_r^{t,a},\bar Y_r^{t,x,a},\bar V_r^{t,x,a})=f(X_r^{t,x,a},I_r^{t,a},\bar Y_r^{t,x,a})$, due to $\bar V_r^{t,x,a}=0$, and $\bar{K}^{t,x,a}$ is a nondecreasing process. Then the second statement readily follows from \cite[Theorem 1.3]{peng99}.
\ep

\vspace{3mm}

The aim is  to obtain  a uniform bound on $(\|Y^{n,t,x,a}\|_{\S^2(t,T)})_n$, which, together with the monotonicity property stated in Lemma \ref{L:Monotone}(i), allows to construct $Y^{t,x,a}$ as the limit of $(Y^{n,t,x,a})_n$. In contrast to \cite[Lemma 3.1]{kharroubi_pham12}, where the compactness of the space of control actions $A$ (which does not hold in the present case, since $A=\R^d$) is exploited to prove the $\S^2$-bound, we utilize a dual (or randomized) representation of $Y^{n,t,x,a}$. To this end, we introduce some additional notations. Let $\mathcal{D}_n$ denote the set of $\mathcal{B}_n$-valued predictable processes $\nu$, where $\mathcal{B}_n$ is the closed ball in $\R^d$ with radius $n$ and centered at the origin, and define $\mathcal{D}= \cup_n \mathcal{D}_n$. For $t\in [0,T]$ and $\nu \in \mathcal{D}$, define the probability measure
$\mathbb{P}^\nu$ on $(\Omega,\Fc_T^t)$ via
\[
 \frac{d\mathbb{P}^\nu}{d\mathbb{P}} = \mathcal{E}\left(\int_t^\cdot \nu_s dB_s\right)_T = \exp\left(\int_t^T \nu_s dB_s - \tfrac{1}{2} \int_t^T |\nu_s|^2 ds\right),
\]
where the Dol\'{e}ans-Dades stochastic exponential on the right-hand side is a martingale due to the definition of $\mathcal{D}$.  In the sequel, we denote by
$\E_s^{t,\nu}$ the condition expectation under $\P^\nu$ given $\Fc_s^t$, for any $s$ $\in$ $[t,T]$, and $\nu$ $\in$ $\mathcal{D}$.

\begin{Remark} \label{remestimXalpha}
{\rm Notice by Girsanov's  theorem that $W$ remains a Brownian motion under $\P^\nu$ for any $\nu\in\mathcal{D}$. Then,
by the same argument as in the derivation of \reff{estimXprho<1} in Lemma \ref{Estimate_X}, we see that when $p_\varrho < 1$,
for every $m\geq 1$, there exists a positive constant $C$, depending only on $m$, $T$, $L_{b,\sigma,\varrho}$,  $M_\varrho$, and $p_\varrho$ such that
\beq \label{estimXalphapro<1}
\E^{t,\nu}_s\big[\sup_{r\in[s,T]}|X_r^{t,x,a}|^m\big] & \leq &
C\bigg(1 + |X_s^{t,x,a}|^m + \int_s^T \E^{t,\nu}_s\big[|I_r^{t,a}|^{\frac{m}{1-p_\varrho}}\big] dr\bigg),
\enq
for all $(t,x,a)\in[0,T]\times\R^d\times \R^d$, $s\in[t,T]$, and $\nu\in\mathcal{D}$.
\epR
}
\end{Remark}

\vspace{1mm}

The next result provides a dual representation of the solution to the penalized BSDE.

\begin{Proposition}
\label{P:DualFormula}
For every $(t,x,a)\in[0,T]\times\R^d\times\R^d$ and $n\in\N$, it holds that
\begin{equation}
\label{DualFormula}
Y_s^{n,t,x,a} \ = \ \underset{\nu\in\mathcal{D}_n}{\textup{ess\,inf}}\,\,\E_s^{t,\nu}\bigg[\int_s^T e^{\int_s^r \gamma_u^n du} f(X_r^{t,x,a},I_r^{t,a},0) dr + e^{\int_s^T \gamma_u^n du} g(X_T^{t,x,a})\bigg],
\end{equation}
for all $t\leq s\leq T$, where the predictable process $\gamma^n\colon\Omega\times[t,T]\rightarrow\R$ is given by
\[
\gamma^n \ = \ \frac{f(X^{t,x,a},I^{t,a},Y^{n,t,x,a}) - f(X^{t,x,a},I^{t,a},0)}{Y^{n,t,x,a}} 1_{\{Y^{n,t,x,a}\neq0\}}.
\]
In particular, $|\gamma^n|$ is  bounded uniformly by the constant $L_F$ in \eqref{Lipschitz_f}.
\end{Proposition}
\textbf{Proof.}
Applying It\^o's formula to $e^{\int_s^\cdot \gamma_u^n du}Y^{n,t,x,a}$ and integrating between $s$ and $T$,  we obtain
\begin{align*}
Y_s^{n,t,x,a} \ &= \ e^{\int_s^T \gamma_u^n du} g(X_T^{t,x,a}) + \int_s^T e^{\int_s^r \gamma_u^n du} f(X_r^{t,x,a},I_r^{t,a},0) dr -  n\int_s^T e^{\int_s^r \gamma_u^n du} |V_r^{n,t,x,a}| dr\\
&\quad \ - \int_s^T e^{\int_s^r \gamma_u^n du} Z_r^{n,t,x,a} dW_r - \int_s^T e^{\int_s^r \gamma_u^n du} V_r^{n,t,x,a} dB_r.
\end{align*}
Take $\nu\in\mathcal{D}_n$, then add and subtract the term $\int_s^T e^{\int_s^r \gamma_u^n du} V_r^{n,t,x,a} \nu_r dr$, so that we obtain
\begin{equation}
\label{BSDE_n_Proof}
\begin{split}
Y_s^{n,t,x,a} \ &= \ e^{\int_s^T \gamma_u^n du} g(X_T^{t,x,a}) + \int_s^T e^{\int_s^r \gamma_u^n du} f(X_r^{t,x,a},I_r^{t,a},0) dr  \\
&\quad \ - \int_s^T e^{\int_s^r \gamma_u^n du} \big(n |V_r^{n,t,x,a}| + V_r^{n,t,x,a}\nu_r \big) dr \\
& \quad \ - \int_s^T e^{\int_s^r \gamma_u^n du} Z_r^{n,t,x,a} dW_r
- \int_s^T e^{\int_s^r \gamma_u^n du} V_r^{n,t,x,a} dB_r^\nu,
\end{split}
\end{equation}
where $B^\nu = -\int_t^\cdot \nu_s \ ds + B$ and $W$ are both $\mathbb{P}^\nu$-Brownian motions by Girsanov's theorem. Note that both local martingales in \eqref{BSDE_n_Proof} are in fact $\mathbb{P}^\nu$-martingales. Indeed, from Bayes formula and the Cauchy-Schwarz inequality, we have
\[
\E_t^{t,\nu}\bigg[\sqrt{\Big\langle\int_t^\cdot Z_r dW_r\bigg\rangle_T}\bigg] \ = \ \E\bigg[L_T^\nu\sqrt{\Big\langle\int_t^\cdot Z_r dW_r\bigg\rangle_T}\bigg]
\ \leq \ \|L_T^\nu\|_{_{\mathbb{L}^2(\Omega,\Fc_T^t,\P)}} \|Z\|_{_{\mathbb H^2(t,T)}},
\]
where $L_T^\nu$ $=$ $d\P^\nu/d\P$ is bounded in ${\mathbb{L}^2(\Omega,\Fc_T^t,\P)}$ since $\nu$ is bounded.
Combining the previous estimate with the Burkholder-Davis-Gundy inequality, we obtain that $\int_t^{\cdot\wedge T} Z_r dW_r$ is of class $(D)$ (see for instance Definition IV.1.6 in \cite{revuz_yor}) under $\P^\nu$, hence it is a $\P^\nu$-martingale (Proposition IV.1.7 in \cite{revuz_yor}). The same argument can be applied to the other local martingale. Projecting both sides of \eqref{BSDE_n_Proof} on $\Fc^t_s$, and using the inequality
$n|v|+ v \cdot a\geq0$,  for any $v\in\R$ and $a\in B_n$, we obtain
\[
Y_s^{n,t,x,a} \ \leq \ \underset{\nu\in\mathcal{D}_n}{\textup{ess\,inf}}\,\,\E_s^{t,\nu}\bigg[\int_s^T e^{\int_s^r \gamma_u^n du} f(X_r^{t,x,a},I_r^{t,a},0) dr + e^{\int_s^T \gamma_u^n du} g(X_T^{t,x,a})\bigg].
\]

To prove the reverse inequality, denote $V^{n;i}$ the $i$-th component of $V^{n,t,x,a}$, for $i=1, \dots, d$. Set $\nu^{*;i}_r = -n V^{n;i}_r /|V^{n}_r| 1_{\{|V^{n}_r \neq 0|\}}$, for $r\in [t,T]$ and $i=1, \dots, d$. Then $\nu^* \in \mathcal{D}_n$, and $n |V^{n,t,x,a}_r| + V^{n,t,x,a}_r \nu^*_r =0$ for any $r\in[t,T]$. It then follows from  \eqref{BSDE_n_Proof} that
\[
Y_s^{n,t,x,a} \ = \ \E_s^{t,\nu^*}\bigg[\int_s^T e^{\int_s^r \gamma_u^n du} f(X_r^{t,x,a},I_r^{t,a},0) dr + e^{\int_s^T \gamma_u^n du} g(X_T^{t,x,a})\bigg],
\]
from which the claim follows.
\ep

\vspace{2mm}

Combining the above dual representation together with the moment estimate for $X$ in Remark \ref{remestimXalpha}, we obtain the following uniform estimate.

\begin{Corollary}
\label{C:Bounds}
There exists a positive constant $C$, independent of $n$, such that
\beq
 && - C \big(1 + |X_s^{t,x,a}|^{q_F\vee p_g}\big) \leq  Y_s^{n,t,x,a}  \nonumber \\
&& \qquad  \leq   \left\{
\begin{array}{lc}
C \big(1 + |X_s^{t,x,a}|^{p_F\vee q_g} + |I_s^{t,a}|^{\frac{p_F\vee q_g}{1-p_\varrho}\vee p'}\big), & \mbox{ if } p_\varrho <1,  \\
C\bigg(1 +  |X_s^{t,x,a}|^{p_F\vee q_g} + |I_s^{t,a}|^{p_F\vee q_g\vee p'} \bigg) e^{C |I_s^{t,a}|} & \mbox{ if } p_\varrho  = 1,
\end{array}
\right.  \label{E:Bounds}
\enq
for all $(t,x,a)\in[0,T]\times\R^d\times\R^d$, $s\in[t,T]$, and $n\in\N$. In particular, $Y^{n,t,x,a}\in\S^2(t,T)$ and
\beq
\|Y^{n,t,x,a}\|_{_{\S^2(t,T)}} & \leq &  \left\{
\begin{array}{lc}
\tilde C \big(1 + |x|^{p_F \vee q_F\vee p_g \vee q_g} + |a|^{\frac{p_F\vee q_g}{1-p_\varrho}\vee p'}\big),  & \mbox{ if } p_\varrho <1,  \\
\begin{array}{l}
\tilde C \big(1 + |x|^{p_F \vee q_F\vee p_g \vee q_g} + |a|^{p_F\vee q_g\vee p'}\big) e^{\tilde C|a|},
\end{array}
&   \mbox{ if } p_\varrho  = 1,
\end{array}
\right.
\label{E:S2}
\enq
for some positive constant $\tilde C$, independent of $n$.
\end{Corollary}
\textbf{Proof.}
{\bf The upper bound in \eqref{E:Bounds}.} From the dual representation formula \eqref{DualFormula}, we have (recall that, when $\nu\equiv0$, $\P^\nu$ coincides with $\P$)
\[
Y_s^{n,t,x,a} \ \leq \ \E_s^t\bigg[\int_s^T e^{\int_s^r \gamma_u^n du} f(X_r^{t,x,a},I_r^{t,a},0) dr + e^{\int_s^T \gamma_u^n du} g(X_T^{t,x,a})\bigg].
\]
Then, recalling the fact that $|\gamma^n|$ is bounded by $L_F$, exploiting the polynomial growth condition of $f$ in \eqref{GrowthCond_f} and of $g$ in Assumption \ref{standing_assumption} (v), we get
\beq \label{interYn}
Y_s^{n,t,x,a} & \leq & C\bigg(1 +  \E_s^t\bigg[ \sup_{r\in[s,T]} |X_r^{t,x,a}|^{p_F\vee q_g} \;  + \;  \sup_{r\in[s,T]} |I_r^{t,a}|^{p'}  \bigg]  \bigg).
\enq
Together with the estimate \reff{estimXprho<1}, and the standard estimate $\E_s^t[\sup_{r\in[s,T]}|I^{t,a}_r|^m] \leq C(1+ |I^{t,a}_s|^m)$ for any $m\geq 1$, this shows
the upper bound in \eqref{E:Bounds}   when $p_\varrho$ $<$ $1$.

When  $p_\varrho$ $=$ $1$, by using the estimate \reff{estimXprho=1}, together with the fact that $ab\leq (a^2 + b^2)/2$ for any $a,b\in \mathbb{R}$, we obtain from \reff{interYn}:
\beq\label{ExpIneq0}
Y_s^{n,t,x,a} & \leq & C\bigg(1 +  |X_s^{t,x,a}|^{p_F\vee q_g} + |I_s^{t,a}|^{p_F\vee q_g\vee p'} \bigg)\sqrt{\E_s^t\big[e^{C \int_s^T |I_r^{t,a}| dr }\big]}.
\enq
For the expectation on the right-hand side, we have
\begin{align}
&\E_s^t\big[e^{C\int_s^T |I_r^{t,a}| dr}\big] \ \leq \ e^{C(T-s) |I_s^{t,a}|} \E_s^t\big[e^{C\int_s^T |B_r - B_s| dr}\big] \ = \ e^{C(T-s) |I_s^{t,a}|} \E\big[e^{C\int_0^{T-s} |B_t| dr}\big] \notag \\
&\leq \ e^{CT|I_s^{t,a}|} \E\big[e^{CT\sup_{0\leq t\leq T}|B_t|}\big] \ = \ e^{CT|I_s^{t,a}|} \E\big[e^{CT\sup_{0\leq t\leq T}\sqrt{|B_t^1|^2 + \cdots + |B_t^d|^2}}\big] \notag \\
&\leq \ e^{CT|I_s^{t,a}|} \E\big[e^{CT\sup_{0\leq t\leq T}(|B_t^1| + \cdots + |B_t^d|)}\big] \ \leq \ e^{CT|I_s^{t,a}|} \E\big[e^{CT(\sup_{0\leq t\leq T}|B_t^1| + \cdots + \sup_{0\leq t\leq T}|B_t^d|)}\big] \notag \\
&= \ e^{CT|I_s^{t,a}|} \E\big[e^{CT\sup_{0\leq t\leq T}|B_t^1|}\big] \cdots \E\big[e^{CT\sup_{0\leq t\leq T}|B_t^d|}\big], \label{ExpIneq}
\end{align}
where the last equality follows from the independence of $B^1,\ldots,B^d$, the components of the $d$-dimensional Brownian motion $B$. Now, notice that
\begin{align}
&\E\big[e^{CT\sup_{0\leq t\leq T}|B_t^1|}\big] \cdots \E\big[e^{CT\sup_{0\leq t\leq T}|B_t^d|}\big] \notag \\
&\leq \ \E\big[e^{CT\sup_{0\leq t\leq T}(B_t^1)^+}e^{CT\sup_{0\leq t\leq T}(B_t^1)^-}\big] \cdots \E\big[e^{CT\sup_{0\leq t\leq T}(B_t^d)^+}e^{CT\sup_{0\leq t\leq T}(B_t^d)^-}\big] \notag \\
&\leq \ \E\big[e^{2CT\sup_{0\leq t\leq T}(B_t^1)^+}\big]^{1/2} \E\big[e^{2CT\sup_{0\leq t\leq T}(B_t^1)^-}\big]^{1/2} \cdots \notag \\
&\quad \ \cdots \E\big[e^{2CT\sup_{0\leq t\leq T}(B_t^d)^+}\big]^{1/2}\E\big[e^{2CT\sup_{0\leq t\leq T}(B_t^d)^-}\big]^{1/2}. \label{ExpIneq2bis}
\end{align}
Using the property that, for every $j=1,\ldots,d$, $(-B_t^j)_{t\geq0}$ is still a Brownian motion, we see that the stochastic processes $((B_t^j)^+)_{t\geq0}=(\max(B_t^j,0))_{t\geq0}$ and $((B_t^j)^-)_{t\geq0}=(\max(-B_t^j,0))_{t\geq0}$ have the same law. As a consequence, the random variables $\sup_{0\leq t\leq T}(B_t^j)^+$ and $\sup_{0\leq t\leq T}(B_t^j)^-$ have the same distribution. Moreover, the distribution of $\sup_{0\leq t\leq T}(B_t^j)^+$ (or, equivalently, of $\sup_{0\leq t\leq T}(B_t^j)^-$) is independent of $j=1,\ldots,d$. Therefore, $\sup_{0\leq t\leq T}(B_t^j)^+$ and $\sup_{0\leq t\leq T}(B_t^j)^-$ have the same distribution as $\sup_{0\leq t\leq T}(B^1_t)^+$. Then, \eqref{ExpIneq2bis} becomes
\[
\E\big[e^{CT\sup_{0\leq t\leq T}|B_t^1|}\big] \cdots \E\big[e^{CT\sup_{0\leq t\leq T}|B_t^d|}\big] \ \leq \ \E\big[e^{2CT\sup_{0\leq t\leq T}(B^1_t)^+}\big]^d.
\]
As $\sup_{0\leq t\leq T}B^1_t\geq0$, $\P$-a.s., (since $B^1_0=0$, $\P$-a.s.), it follows that, $\P$-a.s., $\sup_{0\leq t\leq T}(B^1_t)^+=\sup_{0\leq t\leq T}B^1_t$. In conclusion, \eqref{ExpIneq} becomes
\begin{equation}\label{ExpIneq3bis}
\E_s^t\big[e^{C\int_s^T |I_r^{t,a}| dr}\big] \ \leq \ e^{CT|I_s^{t,a}|}\E\big[e^{2CT\sup_{0\leq t\leq T}B^1_t}\big]^d.
\end{equation}
Using the reflection principle (Proposition III.3.7 in \cite{revuz_yor}), and in particular that $\sup_{0\leq t\leq T}B^1_t$ has the same law as $|B^1_T|$, we obtain from \eqref{ExpIneq3bis}
\begin{align}
&\E_s^t\big[e^{C\int_s^T |I_r^{t,a}| dr}\big] \ \leq \ e^{CT|I_s^{t,a}|} \E\big[e^{2CT|B^1_T|}\big]^d \ \leq \ Ce^{CT|I_s^{t,a}|} \ \leq \ Ce^{C|I_s^{t,a}|}. \label{ExpIneq2}
\end{align}
Combining the previous two estimates \eqref{ExpIneq0} and \eqref{ExpIneq2}, we confirm the upper bound in \eqref{E:Bounds} for $p_\varrho=1$ case.

\vspace{2mm}

{\bf The lower bound in \eqref{E:Bounds}.} From formula \eqref{DualFormula}, the polynomial growth condition of $f$ and $g$, we find that $Y_s^{n,t,x,a}$ is greater than or equal to the following quantity
\[
\underset{\nu\in\mathcal{D}_n}{\textup{ess\,inf}}\,\,\E_s^{t,\nu} \bigg[- M_F \int_s^T e^{\int_s^r \gamma_u^n du} \big(1 + |X_r^{t,x,a}|^{q_F} - \tfrac{1}{q'M_F^{q'}}|I_r^{t,a}|^{q'}\big) dr - m_g e^{\int_s^T \gamma_u^n du} \big(1 + |X_T^{t,x,a}|^{p_g}\big)\bigg].
\]
In the case where $p_\varrho$ $=$ $1$, we have $q_F$ $=$ $p_g$ $=$ $0$, and since $|\gamma^n|$ is bounded by $L_F$, we then obtain the lower bound:
$Y_s^{n,t,x,a}$ $\geq$ $- 2 (M_F T +m_g) e^{L_F T}$.  When $p_\varrho$ $<$ $1$, we use  the estimate \reff{estimXalphapro<1} and  the fact that $|\gamma^n|$ is bounded by $L_F$, to deduce that there exists a positive constant $\bar C$, depending only on $q_F,p_g,T,m_f, q', L_{b,\sigma,\varrho}, L_F$, $M_\varrho$ and $p_\varrho$, such that
\begin{align*}
Y_s^{n,t,x,a} \ &\geq \ - \bar C \big(1 + |X_s^{t,x,a}|^{q_F} + |X_s^{t,x,a}|^{p_g}\big) \\
&\quad \ + \; \underset{\nu\in\mathcal{D}_n}{\textup{ess\,inf}}\,\,\E_s^{t,\nu}\bigg[ \int_s^T \big(\tfrac{M_F^{1-q'}}{q'} e^{-L_F T} |I_r^{t,a}|^{q'} - \bar C |I_r^{t,a}|^{\frac{q_F}{1-p_\varrho}}
- \bar C |I_r^{t,a}|^{\frac{p_g}{1-p_\varrho}}\big) dr \bigg].
\end{align*}
Recalling that $q_F,p_g$ $<$ $(1-p_\varrho)q'$, we see that the function $h(x)$ $=$ $\tfrac{M_F^{1-q'}}{q'} e^{-L_F T}|x|^{q'}$ $-$ $\bar C|x|^{\frac{q_F}{1-p_\varrho}}$ $-$ $\bar C|x|^{\frac{p_g}{1-\varrho}}$,
$x\in\R^d$, is bounded from below on $\R^d$ by a constant independent of $n$. Therefore, we obtain the lower bound.

\vspace{2mm}

{\bf Estimate \eqref{E:S2}.} When $p_\varrho<1$, it follows directly from bounds \eqref{E:Bounds}, Lemma \ref{Estimate_X} with $s=t$,  the standard estimate
$\E[\sup_{s\in[t,T]}|I^{t,a}_s|^m] \leq C(1+ |a|^m)$ for any $m\geq 1$. When
$p_\varrho$ $=$ $1$, \eqref{E:S2} follows from the Cauchy-Schwarz inequality together with $\E[e^{C |I_s^{t,a}|}]$ $\leq$ $ Ce^{C|a|}$ (this latter inequality follows from $\E[e^{C |I_s^{t,a}|}]\leq e^{C|a|}\E[e^{C |B_s-B_t|}]=e^{ C|a|}\E[e^{C |B_{s-t}|}]\leq e^{C|a|}\E[e^{C \sup_{0\leq t\leq T}|B_t|}]$, afterwards we reason as in \eqref{ExpIneq}).
\ep

\vspace{2mm}

The previous uniform norm estimate implies the following uniform norm estimate on $(Z^{n,t,x,a}, V^{n,t,x,a}, K^{n,t,x,a})_n$.

\begin{Corollary}
There exists a positive constant $C$, independent of $n$, such that
\beq
& & \|Z^{n,t,x,a}\|_{_{\mathbb H^2(t,T)}} + \|V^{n,t,x,a}\|_{_{\mathbb H^2(t,T)}} + \|K^{n,t,x,a}\|_{_{\mathbb S^2(t,T)}} \nonumber \\
& &\qquad \leq
\left\{
\begin{array}{lc}
C \big(1 + |x|^{p_F \vee q_F\vee p_g \vee q_g} + |a|^{\frac{p_F\vee q_g}{1-p_\varrho}\vee p'}\big),  & \mbox{ if } p_\varrho <1  \\
\begin{array}{l}
C \big(1 + |x|^{p_F \vee q_F\vee p_g \vee q_g} + |a|^{p_F\vee q_g\vee p'}\big) e^{C|a|},
\end{array}
&   \mbox{ if } p_\varrho  = 1,
\end{array}
\right. \label{EstimateZVK}
\enq
for all $(t,x,a)\in[0,T]\times\R^d\times\R^d$, $n\in\N$.
\end{Corollary}
{\bf Proof.}  By proceeding along the same arguments as in the proof of \cite[Lemma 2.3]{kharroubi_pham12}, we have:
\begin{align}
&\|Z^{n,t,x,a}\|_{_{\mathbb H^2(t,T)}}^2 + \|V^{n,t,x,a}\|_{_{\mathbb H^2(t,T)}}^2 + \|K^{n,t,x,a}\|_{_{\mathbb S^2(t,T)}}^2 \notag  \nonumber \\
& \qquad \leq \ C\bigg(\E\left[\int_t^T |f(X_s^{t,x,a},I_s^{t,a},0)|^2 ds \right] + \sup_{n\in\N}\|Y^{n,t,x,a}\|_{_{\mathbb S^2(t,T)}}^2\bigg), \nonumber
\end{align}
Then, by exploiting the polynomial growth condition of $f$ in \eqref{GrowthCond_f} combined again with Lemma \ref{Estimate_X} for  $s=t$,  the standard estimate
$\E[\sup_{r\in[t,T]}|I^{t,a}_r|^m] \leq C(1+ |a|^m)$ for any $m\geq 1$, and $\E[e^{C |I_s^{t,a}|}]$ $\leq$ $Ce^{C|a|}$ when $p_\varrho$ $=$ $1$,
and using the estimate  \eqref{E:S2} for $Y^{n,t,x,a}$, we obtain the required uniform estimate.
\ep

\vspace{2mm}

Now the previous uniform norm estimates allow us to take limit as $n\rightarrow \infty$.
\begin{Theorem}
\label{T:Y}
Let Assumption \ref{standing_assumption} hold. For every $(t,x,a)\in[0,T]\times\R^d\times\R^d$, there exists a unique maximal solution $(Y^{t,x,a},Z^{t,x,a},V^{t,x,a},K^{t,x,a})\in\mathbb S^2(t,T)\times\mathbb H^2(t,T)\times\mathbb H^2(t,T)\times\mathbb K^2(t,T)$ to \eqref{BSDE}-\eqref{BSDE_Constraint}, such that, as $n\rightarrow\infty$,  $Y_s^{n,t,x,a}\searrow Y_s^{t,x,a}$; $K_s^{n,t,x,a}$ weakly converges to $K_s^{t,x,a}$ in $L^2(\Omega,\Fc_s,\P)$, for all $s\in [t,T]$; and $(Z^{n,t,x,a},V^{n,t,x,a})_n$ weakly converges to $(Z^{t,x,a},V^{t,x,a})$ in $\mathbb H^2(t,T)\times\mathbb H^2(t,T)$. Moreover
\beq
\|Y^{t,x,a}\|_{_{\S^2(t,T)}} & \leq &  \left\{
\begin{array}{lc}
\tilde C \big(1 + |x|^{p_F \vee q_F\vee p_g \vee q_g} + |a|^{\frac{p_F\vee q_g}{1-p_\varrho}\vee p'}\big),  & \mbox{ if } p_\varrho <1  \\
\begin{array}{l}
\tilde C \big(1 + |x|^{p_F \vee q_F\vee p_g \vee q_g} + |a|^{p_F\vee q_g\vee p'}\big) e^{\tilde C|a|},
\end{array}
&   \mbox{ if } p_\varrho  = 1
\end{array}
\right.
\label{E:S2bis}
\enq
where $\tilde C$ is the same constant as in \eqref{E:S2}.
\end{Theorem}

\noindent\textbf{Proof.}
The proof follows the same passages of the proof for \cite[Theorem 2.1]{kharroubi_pham12}, with some small modifications due to the fact that here we adopted a Brownian (rather than Poisson) randomization. For this reason, we just outline main steps of the proof. For every $(t,x,a)\in[0,T]\times\R^d\times\R^d$, from Lemma \ref{L:Monotone}(i) and estimate \eqref{E:S2}, it follows that the sequence $(Y^{n,t,x,a})_n$ converges decreasingly to some $\F^t$-adapted process $Y^{t,x,a}$ satisfying the bound \eqref{E:S2bis}. Next, the monotonic limit theorem (see \cite[Theorem 2.4]{peng99}) implies the weak convergence stated in the theorem and that the limit $(Y^{t,x,a},Z^{t,x,a},V^{t,x,a},K^{t,x,a}) \in \mathbb S^2(t,T)\times\mathbb H^2(t,T)\times\mathbb H^2(t,T)\times\mathbb K^2(t,T)$ satisfies BSDE \eqref{BSDE}. In order to prove the constraint \eqref{BSDE_Constraint}, define the functional $F(V)=\E\int_t^T|V_s|ds$, for $V\in\mathbb H^2(t,T)$. Notice that $\E K_T^{n,t,x,a}=nF(V^{n,t,x,a})$. Hence estimate  \eqref{EstimateZVK} implies  $\sup_n\E|K_T^{n,t,x,a}|^2<\infty$, therefore $\lim_{n\rightarrow \infty}F(V^{n,t,x,a})=0$. Since $F$ is lower semicontinuous with respect to the weak topology of $\mathbb H^2(t,T)$ ($F$ is convex and strongly continuous), it then implies $F(V^{t,x,a})=0$ and \eqref{BSDE_Constraint} holds. Finally, the maximality of $(Y^{t,x,a},Z^{t,x,a},V^{t,x,a},K^{t,x,a})$ follows from a direct application of Lemma \ref{L:Monotone}(ii) and the uniqueness follows from \cite[Proposition 1.6]{peng99}.
\ep

\section{Feynman-Kac representation formula}
\label{S:Feynman-Kac}

\setcounter{equation}{0}
\setcounter{Assumption}{0}
\setcounter{Theorem}{0}
\setcounter{Proposition}{0}
\setcounter{Corollary}{0}
\setcounter{Lemma}{0}
\setcounter{Definition}{0}
\setcounter{Remark}{0}
\setcounter{Example}{0}

The present section is devoted to the proof of Theorem \ref{T:Feynman-Kac}. Assumption \ref{standing_assumption} is in force throughout this section. Let us first recall the definition of (discontinuous) viscosity solution to equation \eqref{VHJ} (or, equivalently, to \eqref{VHJ_0}). Define the \emph{lower semicontinuous (lsc) envelope} $u_*$ and \emph{upper semicontinuous (usc) envelope} $u^*$ of a locally bounded function $u\colon[0,T)\times\R^d\rightarrow\R$ as follows:
\beqs
u_*(t,x) = \liminf_{\substack{(t',x')\rightarrow (t,x) \\ t' < T}} u(t',x') \qquad \text{ and } \qquad u^*(t,x) = \limsup_{\substack{(t',x')\rightarrow (t,x) \\ t' < T}} u(t',x')
\enqs
for all $(t,x) \in [0,T]\times\R^d$.

\begin{Definition}
\label{D:Visc}
\quad
\begin{itemize}
\item[\textup{(i)}] We say that an usc $($resp. lsc$)$ function $u\colon[0,T]\times\R^d\rightarrow\R$ is a viscosity subsolution $($resp. supersolution$)$ to \eqref{VHJ} if
\[
u(T,x) \ \leq \ (resp. \; \geq) \  g(x), \qquad \forall\,x\in\R^d
\]
and
\[
-\frac{\partial\varphi}{\partial t}(t,x) - \inf_{a\in\R^d}\big[\Lc^a \varphi(t,x) + f(x,a,u(t,x))\big] \ \leq \ (resp. \; \geq) \  0
\]
for every $(t,x)\in[0,T)\times\R^d$ and every $\varphi\in C^{1,2}([0,T]\times\R^d)$ satisfying
\[
(u - \varphi)(t,x) \  = \ \max_{[0,T]\times\R^d} (u - \varphi) \quad \big(resp. \ \min_{[0,T]\times\R^d} (u - \varphi)\big).
\]
\item[\textup{(ii)}] We say that a locally bounded function $u\colon[0,T)\times\R^d\rightarrow\R$ is a (discontinuous) viscosity solution to \eqref{VHJ} whenever $u^*$ is a viscosity subsolution and $u_*$ is a viscosity supersolution to \eqref{VHJ}.
\end{itemize}
\end{Definition}

\subsection{Penalized BSDE and corresponding semilinear parabolic PDE}

Let us define, for every $n\in\N$, functions $v_n, v\colon[0,T]\times\R^d\times\R^d\rightarrow\R$ via
\[
v_n(t,x,a) \ = \ Y_t^{n,t,x,a} \quad \text{and} \quad v(t,x,a) \ = \ Y_t^{t,x,a}, \quad \text{ for } \,(t,x,a)\in[0,T]\times\R^d\times\R^d.
\]
Notice that $v_n\searrow v$ pointwise as $n$ goes to infinity. Moreover, it follows from
\eqref{E:S2} and \eqref{E:S2bis} that
\begin{equation}
|v_n(t,x,a)| + |v(t,x,a)| \leq \left\{
\begin{array}{lc}
\tilde C \big(1 + |x|^{p_F \vee q_F\vee p_g \vee q_g} + |a|^{\frac{p_F\vee q_g}{1-p_\varrho}\vee p'}\big),  & \mbox{ if } p_\varrho <1  \\
\begin{array}{l}
\tilde C \big(1 + |x|^{p_F \vee q_F\vee p_g \vee q_g} + |a|^{p_F\vee q_g\vee p'}\big) e^{\tilde C|a|},
\end{array}
&   \mbox{ if } p_\varrho  = 1
\end{array}
\right.
\label{PolGrowthCond_n}
\end{equation}
for all $(t,x,a)\in[0,T]\times\R^d\times\R^d$ and  $n\in \mathbb{N}$.
For each $n\in\N$, let us consider the following semilinear parabolic PDE
\begin{equation}
\label{VHJ_n}
\begin{cases}
\displaystyle - \frac{\partial v_n}{\partial t}(t,x,a) - \tfrac{1}{2}\Delta_a v_n(t,x,a) - \Lc^a v_n(t,x,a) \\
- \, f(x,a,v_n(t,x,a)) + n|D_a v_n(t,x,a)| \ = \ 0, &\quad (t,x,a)\in[0,T)\times\R^d\times\R^d, \\
v_n(T,x,a) \ = \ g(x), &\quad (x,a)\in\R^d\times\R^d,
\end{cases}
\end{equation}
where $\Delta_a$ is the Laplace operator with respect to $a$ and $\Lc^a$ is given in \eqref{E:La}.
Then, we have the following result.

\begin{Proposition}
\label{P:ViscProp_v_n}
The function $v_n$ is a continuous viscosity solution to \eqref{VHJ_n}, i.e., $v_n=(v_n)_*=(v_n)^*$ and properties in Definition \ref{D:Visc} are satisfied where the equation is replaced by \eqref{VHJ_n}.
\end{Proposition}
\textbf{Proof.}
When $p_\varrho<1$, so that $v_n$ satisfies a polynomial growth condition, the result is well-known, see e.g. \cite[Theorem 4.3]{pardoux_peng92}. When $p_\varrho=1$, a change of variable is needed in order to deal with the exponential growth of $v_n$ in $a$. More precisely, define the map $\beta\colon\R\rightarrow\R$ as follows
\[
\beta(x) \ = \
\begin{cases}
e^{x} - 1, &\qquad x\geq\log 2, \\
P(x), &\qquad -\log 2\leq x<\log 2, \\
1 - e^{-x}, &\qquad x<-\log 2,
\end{cases}
\]
where $P(x)=Ax^5+Bx^3+Cx$ (for some real constants $A,B,C$) is a polynomial of fifth degree which realizes a smooth paste of the two other branches of the map $\beta$, so that $\beta\in C^2(\mathbb{R})$ (since $P$ is an odd function, in order to determine $A,B,C$ it is enough to realize a $C^2$-paste with the upper branch of $\beta$) and $\beta$ is increasing, hence invertible. We denote by $\alpha\colon\R\rightarrow\R$ the inverse map of $\beta$, which is given by
\begin{equation}\label{alpha}
\alpha(y) \ = \
\begin{cases}
\log(1 + y), &\qquad y\geq1, \\
P^{-1}(y), &\qquad -1\leq y<1, \\
-\log(1 - y), &\qquad y<-1.
\end{cases}
\end{equation}
For $a=(a_1, \dots, a_d), b=(b_1, \dots, b_d)\in \mathbb{R}^d$, we denote $\alpha(b)= (\alpha(b_1), \dots, \alpha(b_d))\in \mathbb{R}^d$, $\beta(a) = (\beta(a_1), \dots, \beta(a_d))\in \mathbb{R}^d$, and define
\[w_n(t,x,b)=v_n(t,x,\alpha(b)).\] From \eqref{PolGrowthCond_n}, with $p_\varrho=1$, we obtain that there exists a positive constant $c$ such that
\begin{equation}
\label{PolGrowthCond_n_w}
|w_n(t,x,b)| \ \leq \ c \big(1 + |x|^{p_F \vee q_F\vee p_g \vee q_g} + |b|^{p_F\vee q_g\vee p'}\big) |b|^{2\vee \tilde{C}},
\end{equation}
for all $(t,x,b)\in[0,T]\times\R^d\times\R^d$ and  $n\in \mathbb{N}$. Notice that
\begin{align*}
D_a v_n(t,x,a) \ &= \ \bigg(\frac{\partial w_n(t,x,\beta(a))}{\partial b_1} \beta'(a_1),\ldots,\frac{\partial w_n(t,x,\beta(a))}{\partial b_d} \beta'(a_d)\bigg), \\
\Delta_a v_n(t,x,a) \ &= \ \sum_{i=1}^d \frac{\partial^2 w_n(t,x,\beta(a))}{\partial b_i^2} (\beta'(a_i))^2 + \sum_{i=1}^d \frac{\partial w_n(t,x,\beta(a))}{\partial b_i} \beta''(a_i).
\end{align*}
Then, $v_n$ is a continuous viscosity solution to \eqref{VHJ_n} if and only if $w_n$ is a continuous viscosity solution to the following equation:
\begin{equation}
\label{VHJ_n_w}
\begin{cases}
\displaystyle - \frac{\partial w_n}{\partial t}(t,x,b) - \tfrac{1}{2}\Delta_b^\beta w_n(t,x,b) - \Lc^{\alpha(b)} w_n(t,x,b) \\
- \, f(x,\alpha(b),w_n(t,x,b)) + n|D_b^\beta w_n(t,x,b)| \ = \ 0, &\quad (t,x,b)\in[0,T)\times\R^d\times\R^d, \\
w_n(T,x,b) \ = \ g(x), &\quad (x,b)\in\R^d\times\R^d,
\end{cases}
\end{equation}
where
\begin{align*}
D_b^\beta w_n(t,x,b) \ &= \ \bigg(\frac{\partial w_n(t,x,b)}{\partial b_1} \beta'(\alpha(b_1)),\ldots,\frac{\partial w_n(t,x,b)}{\partial b_d} \beta'(\alpha(b_d))\bigg), \\
\Delta_b^\beta w_n(t,x,b) \ &= \ \sum_{i=1}^d \frac{\partial^2 w_n(t,x,b)}{\partial b_i^2} (\beta'(\alpha(b_i)))^2 + \sum_{i=1}^d \frac{\partial w_n(t,x,b)}{\partial b_i} \beta''(\alpha(b_i)).
\end{align*}
Proceeding as in \cite[Theorem 4.3]{pardoux_peng92} we can prove that $w_n$ is a continuous viscosity solution to the previous equation, so the claim follows.
\ep

\subsection{The invariance of the function $v$ with respect to the variable $a$}

We distinguish between the two cases $p_\varrho<1$ and $p_\varrho=1$. Let us consider the two following first-order PDEs:
\begin{align}
\label{Equation_a}
|D_a v(t,x,a)| \ = \ 0, \qquad (t,x,a)\in[0,T)\times\R^d\times\R^d, \\
\label{Equation_a2}
|D_b w(t,x,b)| \ = \ 0, \qquad (t,x,b)\in[0,T)\times\R^d\times\R^d.
\end{align}
Here $w(t,x,b)= v(t,x,\alpha(b))$, where $\alpha(b)$ comes from \eqref{alpha}.
\begin{Lemma}
\label{L:Equation_a}
When $p_\varrho<1$ $($resp. $p_\varrho=1$$)$ the function $v$ $($resp. $w$$)$ is a $($discontinuous$)$ viscosity solution to \eqref{Equation_a} $($resp. \eqref{Equation_a2}$)$.
\end{Lemma}
\textbf{Proof.}
We firstly suppose that $p_\varrho<1$, so that $v_n$ and $v$ satisfy a polynomial growth condition as reported in \eqref{PolGrowthCond_n}. Since the supersolution property clearly holds, let us focus on the subsolution property. To this end, since $v$ is a pointwise limit of the decreasing sequence $(v_n)_n$, it is upper semi-continuous and $v=v^*$. Take $(t,x,a)\in[0,T)\times\R^d\times\R^d$ and $\varphi\in C^{1,2}([0,T]\times\R^d\times\R^d)$ satisfying
$
(v - \varphi)(t,x,a) \ = \ \max_{[0,T]\times\R^d\times\R^d} (v - \varphi).
$
From the polynomial growth property \eqref{PolGrowthCond_n}, we can assume without loss of generality (up to a polynomial perturbation of $\varphi$ for large values of $x$ and $a$) that the previous maximum is strict. Moreover, using standard techniques in viscosity solutions (for similar arguments, see for instance \cite{barles94} or Section 6 in \cite{crandall_ishii_lions92}), we can deduce the existence of a bounded sequence $(t_n,x_n,a_n)_n\in[0,T)\times\R^d\times\R^d$ such that
\begin{align*}
(v_n - \varphi)(t_n,x_n,a_n) \ &= \ \max_{[0,T]\times\R^d\times\R^d} (v_n - \varphi) \quad \text{and} \\
(t_n,x_n,a_n,v_n(t_n,x_n,a_n)) \ &\longrightarrow \ (t,x,a,v(t,x,a)) \quad \text{as } n\rightarrow \infty.
\end{align*}
The viscosity subsolution property of $v_n$ then implies that
\begin{align*}
|D_a\varphi(t_n,x_n,a_n)| \ &\leq \ \frac{1}{n}\bigg(\frac{\partial\varphi}{\partial t}(t_n,x_n,a_n) + \frac{1}{2}\Delta_a\varphi(t_n,x_n,a_n) + \Lc^{a_n}\varphi(t_n,x_n,a_n) \\
&\quad \ + f\big(x_n,a_n,v_n(t_n,x_n,a_n)\big)\bigg).
\end{align*}
Sending $n$ to infinity, and using the continuity of the coefficients $b$, $\sigma$, and $f$, we deduce the claim.

Suppose now that $p_\varrho=1$. Reasoning as in the case $p_\varrho<1$ and using the fact that $\beta'>0$, we can prove that, for every $i=1,\ldots,d$, $w$ is a viscosity solution to the first-order PDE:
\[
\bigg|\frac{\partial w(t,x,b)}{\partial b_i}\bigg| \ = \ 0, \qquad (t,x,b)\in[0,T)\times\R^d\times\R^d.
\]
As a consequence, we deduce that $w$ is a viscosity solution to \eqref{Equation_a2}.
\ep

\vspace{3mm}

We can now state the following result on the independence of $v$ with respect to $a$.

\begin{Proposition}
\label{P:NonDep_a}
The function $v$ satisfies
\[
v(t,x,a) \ = \ v(t,x,a'), \quad \text{ for all } (t,x)\in[0,T]\times\R^d \text{ and } a,a'\in\R^d.
\]
\end{Proposition}
\textbf{Proof.}
Let us firstly consider the case $p_\varrho<1$. Then, this result can be proved proceeding as in \cite[Lemma 3.4 and Proposition 3.2]{kharroubi_pham12}, with only two small modifications. First, in \cite{kharroubi_pham12} the equation $-|D_av(t,x,a)|=0$ is considered rather than \eqref{Equation_a} (this is due to the fact that the function $v$ therein is the \emph{increasing} limit of $(v_n)_n$, since a Hamilton-Jacobi-Bellman equation with a ``$\sup$'' operator, instead of ``$\inf$'', is studied there). However, the proofs of \cite[Lemma 3.4 and Proposition 3.2]{kharroubi_pham12} are not affected (apart from obvious modifications by the presence of the minus sign). Second, in \cite{kharroubi_pham12} the variable $a$ belongs to an open, bounded, and connected subset $\mathring A$ of some Euclidean space $\R^q$, rather than to the entire space $\R^d$ as here. To link to \cite{kharroubi_pham12}, let us consider the open ball $\mathring{\mathcal{B}}_r\subset\R^d$ for any $r>0$. Lemma \ref{L:Equation_a} implies that $v$ is a viscosity solution to equation \eqref{Equation_a} on $[0,T)\times\R^d\times \mathring{\mathcal{B}}_r$, for every $r>0$. At this point we can apply  \cite[Lemma 3.4 and Proposition 3.2]{kharroubi_pham12} to deduce that $v$ satisfies
\[
v(t,x,a) \ = \ v(t,x,a'), \quad \text{ for all } (t,x)\in[0,T]\times\R^d \text{ and } a,a'\in\mathcal{B}_r.
\]
As a result, the claim follows since $r$ is arbitrarily chosen.

Suppose now that $p_\varrho=1$. Then, proceeding along the same lines as in the case $p_\varrho<1$, we conclude that $w$ does not depend on $b$. Since $v(t,x,a)=w(t,x,\beta(a))$, this implies that $v$ does not depend on $a$.
\ep

\subsection{From the BSDE with diffusion constraint to the viscous HJ equation}

Following Proposition \ref{P:NonDep_a}, we define the function $u$, as in \eqref{Feynman-Kac}, which satisfies (due to \eqref{PolGrowthCond_n} with $a=0$)
\[
|u(t,x)| \ \leq \ \tilde C \big(1 + |x|^{p_F \vee q_F\vee p_g \vee q_g}\big), \quad  \text{ for all } (t,x)\in[0,T]\times\R^d.
\]
Let us now prove that $u$ is a viscosity solution to equation \eqref{VHJ} (equivalently, to equation \eqref{VHJ_0}), which, together with Proposition \ref{P:NonDep_a}, concludes the proof of Theorem \ref{T:Feynman-Kac}.

\vspace{3mm}

\noindent\textbf{Proof of Theorem \ref{T:Feynman-Kac}.}
In order to prove that $u$ is a viscosity solution to equation \eqref{VHJ}, we begin by noticing that, as a direct consequence of Proposition \ref{P:ViscProp_v_n}, $v_n$ is a viscosity solution to \eqref{VHJ_n} on $[0,T]\times\R^d\times\mathcal{B}_r$, for any $r>0$. Therefore, as a decreasing limit of $(v_n)_n$, $v$, hence  $u$, is a viscosity solution to the following equation:
\begin{equation}
\label{VHJ_r}
\begin{cases}
\displaystyle - \frac{\partial u}{\partial t}(t,x) - \inf_{a\in\mathcal{B}_r}\big[\Lc^a u(t,x) + f(x,a,u(t,x))\big] \ = \ 0, &\quad (t,x)\in[0,T)\times\R^d, \\
u(T,x) \ = \ g(x), &\quad x\in\R^d.
\end{cases}
\end{equation}
The proof of this result can be done proceeding as in  \cite[Section 3.4]{kharroubi_pham12}, with some small modifications, due to the Poisson, rather than Brownian, randomization used in \cite{kharroubi_pham12}, which affects in part the form of the penalized PDE \eqref{VHJ_n}. From the definition of viscosity solution to \eqref{VHJ_r} we have, for any $(t,x)\in[0,T)\times\R^d$, $r>0$, and any $\varphi\in C^{1,2}([0,T]\times\R^d)$,
\[
-\frac{\partial\varphi}{\partial t}(t,x) - \inf_{a\in B_r}\big[\Lc^a \varphi(t,x) + f(x,a,u(t,x))\big] \ \leq \ (\text{resp.} \; \geq) \  0
\]
whenever
$
(u - \varphi)(t,x) \  = \ \max_{[0,T]\times\R^d} (u - \varphi)$  (resp. $ \min_{[0,T]\times\R^d} (u - \varphi)$).
Noticing that the choice of the test function $\varphi$ is independent of $r$, we deduce
\[
-\frac{\partial\varphi}{\partial t}(t,x) - \inf_{a\in\R^d}\big[\Lc^a \varphi(t,x) + f(x,a,u(t,x))\big] \ \leq \ (\text{resp.} \; \geq) \  0,
\]
which implies that $u$ is a viscosity solution to equation \eqref{VHJ}. Finally, the terminal condition is also satisfied thanks to the terminal condition in \eqref{VHJ_n} and the fact that $v$ is a monotone limit of $(v_n)_n$.
\ep

\small
\bibliographystyle{plain}
\bibliography{biblio}

\end{document}